\documentclass[11pt]{article}

\usepackage{amsfonts,amsmath,amssymb,theorem}
\usepackage[T1]{fontenc}
\usepackage[latin1]{inputenc}
\usepackage{graphicx,graphics}
\usepackage{color}
\usepackage{indentfirst}
\usepackage{geometry}%pour les marges
\global
\global
\theoremstyle{break} \theorembodyfont{\it}

\newtheorem{thm}{\textbf{Theorem}} \theorembodyfont{\sl}
\newtheorem{definition}{\textbf{Definition}}
\newtheorem{proposition}{\textbf{Proposition}}

\newtheorem{lemme}{\textbf{Lemma}}
\newtheorem{rem}{Remark}
\theoremheaderfont{\sc}

\def\N{\mathbb{N}}
\def\Z{\mathbb{Z}}
\def\R{\mathbb{R}}
\def\C{\mathbb{C}}
\def\H{\mathbb{H}}
\def\S{\mathbb{S}}

\def\Im{{\rm{Im }}}
\def\Vol{{\rm{Vol }}}
\def\dvol{{\rm{dvol }}}
\def\Ker{{\rm{Ker}\,}}

\def\Res{{\rm{Res}}}
\def\Ran{{\rm{Ran}}}
\def\Tr{{\rm{Tr}}}
\def\ad{{\rm{ad}}}
\def\det{{\rm{det}}}

\def\Re{{\rm{Re }}}
\def\supp{\mathrm{supp}}

\def\min{\mathrm{Min}}
\def\spec{\mathrm{Spec}}
\def\dim{\mathrm{dim}}
\def\sh{\mathrm{sh}}
\def\ch{\mathrm{ch}}
\def\d{\mathrm{d}}

\begin{document}

\parindent=10pt

\vskip 1cm
\begin{center}
{\textbf {\Large{Isoresonant complex-valued potentials and symmetries.}}}
\end{center}

\begin{center}
{ Aymeric AUTIN
\footnote{Laboratoire de mathématiques Jean Leray (UMR 6629). Université de Nantes, 2 rue de la Houssinière - BP 92208
44322 Nantes Cedex 3, France. \textit{Email} : aymeric.autin@univ-nantes.fr.
}}
\end{center}

\;\;\;

\begin{abstract}
Let $X$ be a connected Riemannian manifold such that the resolvent of the free Laplacian $(\Delta-z)^{-1}, z\in\C\setminus\R^{+},$ has a
meromorphic continuation through $\R^{+}$. The poles of this continuation are called resonances. When $X$ has some symmetries, we
construct complex-valued potentials, $V$, such that the resolvent of $\Delta+V$, which has also a meromorphic continuation, has the same
resonances with multiplicities as the free Laplacian. 

\bigskip

 
\end{abstract}

\section{Introduction and statement of the results}

Let $(X,g)$ be a connected Riemannian manifold with dimension $n\geq 2$. On $X$ we have the free non-negative Laplacian, $\Delta$, acting on functions with domain $H^{2}(X)$, whose spectrum is included
in $\R^{+}$. So, for $z\in\C\setminus\R^{+}$, the resolvent $R_{0}(z):=(\Delta-z)^{-1}$ of the Laplacian is a bounded operator from $L^{2}(X)$ to $H^{2}(X)$. 
We will assume that, the resolvent has a meromorphic continuation through $\R^{+}$ in a domain of $\C$, $D^{+}$. For example, this holds for Euclidean spaces, asymptotically hyperbolic manifolds and manifolds with asymptotically cylindrical ends.

We call \textit{resonance} of $\Delta$ a pole of $R_{0}$ in $D^{+}$, and we write $\Res(\Delta)$, the set of these poles. If $z_{0}\in \Res(\Delta)$,
then, in a neighbourhood of $z_{0}$ in $D^+$, we have a finite Laurent expansion :
\[R_{0}(z)=\sum_{i=1}^{p}(z-z_{0})^{-i}S_{i} + H(z),\]
where $S_{i}$ has a finite rank and $H$ is holomorphic. $p$ is the
\textit{order} of the resonance. We call \textit{multiplicity} of $z_{0}$ the dimension of the resonant space which is the range of $S_{1}$. See \cite{A}.

If we perturb the Laplacian with a potential $V$ and if $V$ is sufficiently decreasing at infinity on $X$, for example compactly supported, then the resolvent of $\Delta+V$, $(\Delta+V-z)^{-1}$ can also be continued meromorphically to $D^{+}$. Then, we can introduce the resonances of $(\Delta+V)$, we write their set $\Res(\Delta+V)$.

For such a $V$, sufficiently decreasing, we have the equality for the essential spectrum, $\sigma_{ess}(\Delta+V)=\sigma_{ess}(\Delta)$, because $V$ is then relatively compact with respect to $\Delta$. So we can wonder how these potentials modify resonances. We reach the main question of this work~:\\

\noindent{\bfseries Do there exist potentials $\mathbf{V}$ such that $\mathbf{\Res(\Delta+V)=\Res(\Delta)}$ ?} \\

We will construct such potentials and we call them \textit{isoresonant}. In term of inverse problem, we can't detect their presence only with the observation of the set of resonances.

Our potentials will be complex-valued and it is crucial. For example it is known that, in $\R^{n}$, $n\geq2$ and even or $n=3$, nontrivial, real valued, smooth and compactly supported potentials create an infinite number of resonances. See \cite{MR} \cite{SZ} \cite{C0} \cite{S}.

We have been inspired by the work of Christiansen in \cite{C2} and \cite{C3}. She constructs in Euclidean spaces $\R^{n}$ ($n\geq 2$) isoresonant complex potentials, i.e. in this case : $\Res(\Delta+V)=\Res(\Delta)=\emptyset$. She uses an action of $\S^{1}$ on $\R^{n}$. I generalize this construction to manifolds which have an isometric action of $\S^{1}$, and I use other symmetries as $(\S^{1})^{m}$ and $SO(n)$. On these manifolds, the free Laplacian already has some resonances, so there is more work to prove the isoresonance of the potentials. To compare ; in the Euclidean space, it is sufficient to prove $\Res(\Delta+V)\subset\Res(\Delta)$ because $\Res(\Delta)=\emptyset$.

We are going to describe the method for the construction and the statements of the results. We assume that $(X,g)$ has an isometric action of $\S^{1}$. This action induces an unitary representation of $\S^{1}$ on $L^{2}(X)$ : 
\begin{eqnarray*}
S^{1}&\longrightarrow& U(L^{2}(X)) \\
e^{i\theta}&\longrightarrow& f\to (x\to f(e^{-i\theta}.x)).
\end{eqnarray*}
Then we can decompose $L^{2}(X)$ according to isotypical components :
\[L^{2}(X)=\overline{\bigoplus_{j\in\Z}^{\perp}L^{2}_{j}(X)},\]
with, for all $j\in\Z$,
\[L^{2}_{j}(X):=\{f\in L^{2}(X) \ ; \ \forall \theta\in[0,2\pi],\ \forall x\in X, \ f(e^{-i\theta}.x)=e^{ij\theta}f(x)\},\]
is the space of $\S^{1}$ homogeneous functions of weight $j$.\\
\indent We take for our isoresonant potentials sums of $\S^{1}$ homogeneous functions with weights of the same sign. Such functions create a shift on the isotypical components of $L^{2}(X)$ : if $V\in L^{\infty}(X)\cap L^{2}_{m}(X)$ and $f\in L^{2}_{j}(X)$ then $Vf\in L^{2}_{j+m}(X)$. On the contrary, the Laplacian stabilizes these isotypical components. Thanks to this shift we will prove the inclusion $\Res(\Delta+V)\subset\Res(\Delta)$, first for truncated $V$, and after for all $V$ thanks to a characterization of resonances as zeros of regularized determinants.\\
\indent On the way, we have to estimate, for all compact $K$, the lower bound of the spectrum of the Dirichlet Laplacian acting on $\S^{1}$ homogeneous functions of weight $j$ supported in $K$ (we denote this space $L^{2}_{j}(K)$). This is an interesting result on its own :

\begin{proposition}
\label{prop0}
Let $K$ be a compact manifold with boundary, having an action of $\S^{1}$ and a metric $g$ such that $\S^{1}$ acts by isometries on $(K,g)$ and $g$ has a product form in a neighborhood of the boundary of $K$. Then there exist strictly positive constants, $C_{1}(K)$ and $C_{2}(K)$, such that, for all $j\in\Z$, we have :
\[C_{1} j^{2}\leq\min \; \spec \; \Delta_{L^{2}_{j}(K)}\leq C_{2} (1+j^{2}).\]  
\end{proposition}

For the other inclusion, $\Res(\Delta)\subset\Res(\Delta+V)$, we use the Agmon's perturbation theory of resonances, developed in \cite{A}. Thanks to this theory we can view resonances as eigenvalues of auxiliary operators and so we can use the Kato's theory in order to study their perturbations.\\
\indent Finally we get the following result, given here in restricted cases for simplicity :

\begin{thm}
\label{thm0}
On the Euclidean space $\R^{n}$ or the hyperbolic space $\H^{n}$, let be the potential
\[V=\sum_{m=1}^{M}V_{m},\]
where $V_{m}\in L^{\infty}(X)$ is compactly supported and $\S^{1}$ homogeneous with weight $m$.

Then, in $\C$, we have $\Res(\Delta+V)=\Res(\Delta)$ with the same multiplicities.
\end{thm}

See the theorem \ref{thm : circular} for the general case with a more general manifold, an infinite sum for $V$, and $V$ not compactly supported.

\begin{rem}
\label{remV0}
Instead of the free Laplacian we can perturb $\Delta+V_{0}$ with $V_{0}$ a real, compactly supported and $\S^{1}$ invariant potential and the result becomes $\Res(\Delta+V_{0}+V)=\Res(\Delta+V_{0})$ with the same multiplicities. We can imagine the perturbation of other operators which respect the decomposition of $L^{2}(X)$ according to the isotypical components.
\end{rem}

The construction of isoresonant potentials using the action of $(\S^{1})^{m}$ is essentially the same as in the case $\S^{1}$ so we don't describe it in this article, but it can be found in \cite{AA}. On the contrary, if we look at the action of $SO(n)$ ($n\geq 3$), as this group is not commutative, we don't have any simple description of the isotypical components. Then we add an hypothesis and assume that we can write 
\[L^{2}(X)=\bigoplus_{k\in\N}L^{2}(\R^{+})\otimes H^{k},\]
where $H^{k}= Ker(\Delta_{S^{n-1}}-k(k+n-2))$, $k\in\N$, is the eigenspace of the Laplacian on the sphere $\S^{n-1}$. Like in the case $\S^{1}$, we are going to construct some $V$ which induces a shift in this decomposition of $L^{2}(X)$. This time $V$ is a sum of highest weight vectors of the representations $H^{k}$ of the complexification of the Lie algebra $\mathfrak{so}_{n}$. Moreover, for the action of $SO(n)$ we don't need to use the proposition \ref{prop0}, which simplifies the proof of the isoresonance. Here I have been inspired by the construction of isospectral potentials by Guillemin and Uribe in \cite{GU}.

These potentials don't modify the set of resonances of the free Laplacian and their multiplicity. We can wonder if, with more information, we would be able to detect them. On this way, I prove that , on $\H^{2}$, there exist some potentials among the family of isoresonant potentials which modify the order of the resonances. On $\H^{2}$, resonances of the free Laplacian are, up to a change of spectral parameter, the negative integers with order $1$. Taking for the hyperbolic plan the model $\R^{+}\times \S^{1}$ with coordinates $(r,\theta)$ and metric 
$g=\d r^{2}+\sh(r)^{2}\d\theta^{2}$, we have

\begin{proposition}
On the hyperbolic plan $\H^{2}$, let $k$ be a strictly positive integer. There exists a potential $V\in\mathcal{F}:=\{V_{m}(r)e^{im\theta} ; m\in\Z\setminus\{0\}, V_{m}\in L^{\infty}_{c}(\R^{+})\}$ such that $-k$ is a resonance of $\Delta+V$ with an order strictly bigger than $1$.
\end{proposition}

In the last part of this article, we construct isoresonant potentials using the $\S^{1}$ action on another example : the catenoid, i.e. $(\R\times \S^{1},\d r^{2}+(r^{2}+a^{2})\d\alpha^{2})$ with $(r,e^{i\alpha})\in\R\times\S^{1}$ and $a\in\R$. We treat this example separately because we can't use the Agmon's theory for defining the continuation of the resolvent and for perturbations of resonances ; instead we have to use a complex scaling method, following \cite{WZ}.

\section{Framework and conditions}

We take a cover $f:\Sigma\rightarrow\Omega$, where $\Omega$ 
is an open set of $\C$, and an unbounded domain $D\subset\Sigma$ such that $f(D)\subset\C\setminus\R^{+}$. We note $R_{0}(\lambda):=(\Delta-f(\lambda))^{-1}$ which is first defined holomorphic in $D$ with values in $\mathcal{L}(L^{2}(X))$ (the space of bounded operators from $L^{2}(X)$ to itself). Let two Banach spaces $B_{0}$ and $B_{1}$ be such that 
\[B_{0}\overset{J_{0}}{\hookrightarrow} L^{2}(X)\overset{J}{\hookrightarrow} B_{1},\]
where $J_0$ and $J$ are continuous injections, $J_{0}(B_0)$ is dense in $L^{2}(X)$ and $J(L^{2}(X))$ is dense in $B_1$.
We note, for $\lambda\in D$,
\[\widetilde{R}_{0}(\lambda)=JR_{0}(\lambda)J_{0}.\]
$\widetilde{R}_{0}$ is holomorphic on $D$ with values in $\mathcal{L}(B_{0},B_{1})$.

\vskip 0.5cm

Our first assumption is the following :

\vskip 0.3cm

\textbf{\underline{Condition $A$}} : $\widetilde{R}_{0}$ has a meromorphic continuation with finite rank poles (we will say \textit{finite-meromorphic}) from $D$ to $D^{+}$ a domain of $\Sigma$.

\vskip 0.3cm

In order that Condition A not hold trivially we assume that $f(D^{+})$ intersects the essential spectrum of $\Delta$.

Let us give some examples :

\begin{itemize}
\item $X$ is $\R^{n}$ with the Euclidean metric, 
	\begin{itemize}
	\item[-] If $n$ is odd then we take $\Sigma=\C$ and $R_{0}(\lambda):=(\Delta-\lambda^{2})^{-1}$ is first defined in $D=\{\lambda\in\C\ ;\ \Im\lambda>0\}$ with values in $\mathcal{L}(L^{2}(\R^{n}))$ and has, for all $N>0$, an holomorphic continuation in $D^{+}_{N}=\{\lambda\in\C\ ;\ \mid\!\Im\lambda\!\mid<N\}$ with values in $\mathcal{L}(e^{-N<z>}L^{2}(\R^{n}),e^{N<z>}L^{2}(\R^{n}))$, where $<\!z\!>=(1+\mid\!z\!\mid^{2})^{\frac{1}{2}}$. See \cite{MR} and \cite{SZ}.
	\item[-] If $n$ is even then we take for $\Sigma$ the logarithmic cover of $\C\setminus\{0\}$, and $R_{0}(\lambda):=\\ (\Delta-e^{2\lambda})^{-1}$ is first defined in $D=\{\lambda\in\C\ ;\ 0<\Im\lambda<\pi\}$ with values in $\mathcal{L}(L^{2}(\R^{n}))$ and has, for all $N>0$, an holomorphic continuation in $D^{+}_{N}=\{\lambda\in~\C\ ;\ \mid~\!\Im(e^{\lambda})\!\mid<~N\}$ with values in $\mathcal{L}(e^{-N<z>}L^{2}(X),e^{N<z>}L^{2}(X))$. See \cite{MR}.
	\end{itemize}
	
\item $X$ is an asymptotically hyperbolic manifold. We begin with the definition of such  a manifold. Let $\overline{X}=X\cup\partial\overline{X}$, a smooth compact manifold of dimension $n$ with boundary $\partial\overline{X}$ and $\rho_{0}$ a boundary-defining function that is a smooth function on $\overline{X}$ such that 
\[\rho_{0}\geq 0,  \;\;  \partial\overline{X}=\{m\in\overline{X} \ ; \ \rho_{0}(m)=0\},  \;\;  \d\rho_{0}\mid_{\partial\overline{X}}\neq 0.\]
We say that a smooth metric $g$ on $X$ is asymptotically hyperbolic if $\rho_{0}^{2}g$ continues as a smooth metric on $\overline{X}$ and $\mid\!\d\rho_{0}\!\mid_{\rho_{0}^{2}g}=1$ on $\partial\overline{X}$. Thanks to this condition, the sectional curvature of $g$ tends to $-1$ at the boundary, and there exists a function $\rho$ defining the boundary, a collar neighborhood of the boundary, $U_{\rho}:=[0,\varepsilon)\times\partial\overline{X}$, and a family $h(\rho), \rho\in[0,\varepsilon)$, of smooth metrics on $\partial\overline{X}$ such that
\begin{equation}
\label{metric}
g=\frac{\d\rho^{2}+h(\rho)}{\rho^{2}} \ \ \ \textrm{on $U_{\rho}$.}
\end{equation}
For example the hyperbolic space $\H^{n}$ and its convex co-compact quotients are asymptotically hyperbolic.

We take $\Sigma=\C$ and $R_{0}(\lambda):=(\Delta-\lambda(n-1-\lambda))^{-1}$ is first defined and meromorphic in $D=\{\lambda\in\C\ ;\ \Re\lambda>\frac{n-1}{2}\}$ with values in $\mathcal{L}(L^{2}(X))$ and $\lambda$ is one of its poles if and only if $\lambda(n-1-\lambda)\in\sigma_{d}(\Delta)$, and it is of finite rank. \\
\indent Mazzeo and Melrose (\cite{MM}) and after Guillarmou (\cite{GC}) have proved that $R_{0}$ has a finite-meromorphic continuation in $\C\setminus(\frac{n}{2}-\N)$ and in all $\C$ if and only if the metric $g$ is even. The metric $g$ is \textit{even} if the family $h(\rho)$ defined in (\ref{metric}) has a Taylor's series in $\rho=0$ only with even powers of $\rho$ (it does not depend on the choice of $\rho$). More precisely, for all $N\geq 0$, $R_{0}$ has a finite-meromorphic continuation on $D^{+}_{N}:=\{\lambda\in\C\ ;\ \Re\lambda>\frac{n-1}{2}-N\}$ if $g$ is even, and otherwise on $D^{+}_{N}\setminus(\frac{n}{2}-\N)$, with values in $\mathcal{L}(\rho^{N}L^{2}(X),\rho^{-N}L^{2}(X))$.

\item $X$ is a Riemannian manifold with asymptotically cylindrical ends. Like in the previous case, let $\overline{X}=X\cup\partial\overline{X}$ be a smooth compact manifold of dimension $n$ with boundary. We say that a smooth metric $g$ on $X$ is a metric with asymptotically cylindrical ends if there exists a function $\rho$ defining the boundary, a collar neighborhood of the boundary, $U_{\rho}:=[0,\varepsilon)\times\partial\overline{X}$, and a family $h(\rho), \rho\in[0,\varepsilon)$, of smooth metrics on $\partial\overline{X}$ such that
\begin{equation}
g=\frac{\d\rho^{2}}{\rho^{2}}+h(\rho) \ \ \ \textrm{on $U_{\rho}$.}
\end{equation}
Let $\Delta_{\partial\overline{X}}$ be the Laplacian on the compact manifold $\partial\overline{X}$ and $0=\sigma_{1}^{2}<\sigma^{2}_{2}<\ldots$ its spectrum. Then the spectrum of the Laplacian on $X$, $\Delta$, is discrete with finite multiplicities outside $[0,+\infty)$ and, for all $j>0$, $[\sigma_{j},\sigma_{j+1})$ is continuous spectrum with a multiplicity equal to the sum of the multiplicities of $\{\sigma_{1},\ldots,\sigma_{j}\}$ as eigenvalues of $\Delta_{\partial\overline{X}}$ and there may be embedded eigenvalues with finite multiplicity.

Melrose, in \cite{MV}, proves the continuation of the resolvent of the free Laplacian on the Riemannian surface $\Sigma$ which is the surface such that all the functions $r_{j}(\lambda):=(\lambda-\sigma_{j}^{2})^{\frac{1}{2}}$ are holomorphic on it. This surface is ramified at points $\lambda=\sigma_{j}^{2}$. $R_{0}(\lambda)=(\Delta-\lambda)^{-1}$ is first defined in $D=\{\lambda\in\Sigma\ ;\ \forall \, j  \;  \Im(r_{j}(\lambda))>0\}$ with values in $\mathcal{L}(L^{2}(X))$ and, for all $N\geq 0$, it has a finite-meromorphic continuation in a domain $D^{+}_{N}$ with values in $\mathcal{L}(\rho^{N}L^{2}(X),\rho^{-N}L^{2}(X))$. 

\item $X$ is a rank-one symmetric space of the noncompact type. In this case, Hilgert and Pasquale prove, in \cite{HP}, the finite-meromorphic continuation of the resolvent of the free Laplacian. 
\end{itemize}

\vskip 0.5cm

In order to treat all these examples with a single notation, we reformulate the Condition A, with $N>0$, and $\rho=e^{-<z>}$ for the Euclidean case and a boundary-defining function on $\overline{X}$ for the other examples, by

\vskip 0.3cm

\textbf{\underline{Condition $A_{N,\rho}$}} : $\widetilde{R}_{0}$ has a finite-meromorphic continuation from $D$ to $D^{+}_{N}$ an unbounded domain of $\Sigma$, with values in $\mathcal{L}(\rho^{N}L^{2}(X),\rho^{-N}L^{2}(X))$.

\vskip 0.3cm

Agmon shows in \cite{A} that notions of resonances, multiplicity and order do not depend of the weight $\rho^{N}$ chosen.

\vskip 0.5cm

In order to have the finite-meromorphic continuation to $D_{N}^{+}$ of the resolvent of $\Delta+V$, $(\Delta+V-z)^{-1}$, $z\in\C\setminus\spec(\Delta+V)$, we introduce a condition for $V$ :

\vskip 0.3cm

\textbf{\underline{Condition $B_{N,\rho}$}} : $\widetilde{R}_{V}(\lambda):=J(\Delta+V-f(\lambda))^{-1}J_{0}$ with values in $\mathcal{L}(\rho^{N}L^{2}(X),\rho^{-N}L^{2}(X))$ has a finite-meromorphic continuation from $D$ to $D^{+}_{N}$ and $\rho^{-2N}V$ is bounded on $X$.

\vskip 0.3cm

\begin{rem}
With the hypothesis $\rho^{-2N}V$ bounded, we will be able to apply the Agmon's perturbation theory of resonances \cite{A}. 
\end{rem}

\begin{rem}
If $V$ is compactly supported or if $V$ is smooth on $\overline{X}$ and vanishes to all orders in $\rho$ at the boundary then $V$ verifies Condition $B_{N,\rho}$ for all $N$. In these cases, the resolvent of $\Delta +V$ has a finite-meromorphic continuation in all $\Sigma$.
\end{rem}

\section{Circular symmetries}
\subsection{Statement of the result and examples}

Let $(X,g)$ be a connected Riemannian manifold with a $\S^{1}$ action by isometries. We will need the following condition.

\vskip 0.3cm

\textbf{\underline{Condition $C$}} : For all compact $K\subset X$ there exists a compact manifold with boundary $\widetilde{K}$ which is diffeomorphic to a compact of $X$ containing $K$, and which has a isometric $\S^{1}$ action, with a smooth metric $\widetilde{g}$ such that $\widetilde{g}_{\mid\!K}=g_{\mid\!K}$, and $\widetilde{g}$ is a product metric $\d\delta^{2}+\widetilde{h}_{\partial\widetilde{K}}$ in a neighbourhood of the boundary $\partial\widetilde{K}$ of $\widetilde{K}$, with $\delta$ a $\S^{1}$ invariant function defining $\partial\widetilde{K}$ and $\widetilde{h}$ is independent of $\delta$.

\vskip 0.3cm
 
We can now give the main result of this part :

\begin{thm}
\label{thm : circular}
Let $(X,g)$ be a connected Riemannian manifold with an action of $\S^{1}$ by isometries verifying the Condition C and the Condition $A_{N,\rho}$ for some $N>0$ and some $\S^{1}$ invariant function $\rho$.

Let $V$ be the potential
\[V=\sum_{m=1}^{+\infty}V_{m},\]
where $V_{m}\in L^{\infty}(X)$ is $S^{1}$ homogeneous of weight $m$, with $\underset{m=1}{\overset{+\infty}{\sum}}\parallel\!V_{m}\!\parallel_{\infty}<+\infty$. If $V$ verifies Condition $B_{N,\rho}$, and for all $\lambda\in D_{N}^{+}\setminus\Res(\Delta)$, $\rho^{-(N+1)}V\widetilde{R}_{0}(\lambda)\rho^{N}$ is in a Schatten class $\mathcal{S}_{q}$, $q\in\N\setminus\{0\}$,

\vskip 0.5cm

\noindent then, on  $D_{N}^{+}$, $\Res(\Delta +V)=\Res(\Delta)$ with the same multiplicities.
\end{thm}

We will recall the definition of the Schatten classes in the definition \ref{def : schatten} in section $3.3.2$.

\begin{rem}
The last assumption, $\rho^{-(N+1)}V\widetilde{R}_{0}(\lambda)\rho^{N}\in\mathcal{S}_{q}$, is technical and will allow us to use regularized determinants. 
If $V$ is compactly supported then it holds with $q>\frac{\dim X}{2}$ and for any $N$. In the Euclidean space $\R^{n}$, if $V$ is bounded and super-exponentially decaying we still have the assumption for all $N$. On an asymptotically hyperbolic manifold $X$, if $V$ is smooth on $\overline{X}$ and vanishes to all orders in $\rho$ at the boundary then $V$ verifies the assumption for all $N$.
\end{rem}

\vskip 0.3cm

Let us describe the $\S^{1}$ action, the Condition $C$ and the potentials for the examples of section $2$ verifying Condition $A$ :

 \begin{itemize}
\item Let $\R^{n}=(\R^{2})^{k}\times\R^{n-2k}$ be the Euclidean space with the following $\S^{1}$ action

\[\bigoplus_{i=1}^{k}R(p_{i}\theta)\oplus Id_{\R^{n-2k}},\]
where $\theta\in[0,2\pi)$, $(p_{1},\ldots,p_{k})\in(\Z\setminus\{0\})^{k}$, and $R(\phi)$ is the rotation of angle $\phi$ on $\R^{2}$. For the Condition $C$, we can remark that every compact $K$ can be included in a ball $B(0,R)$ and we can take for $\widetilde{K}$ a bigger ball $B(0,\widetilde{R})$ with $\widetilde{R}>R$ with the following metric in polar coordinates :
\[\widetilde{g}=dr^{2}+f(r)d\omega^{2}, \;\; (r,\omega)\in\R^{+}\times\S^{n-1}, \]
where $d\omega^{2}$ is the metric on the $(n-1)$-sphere in $\R^{n}$ and $f$ is smooth on $[0,\widetilde{R}]$, constant near $\widetilde{R}$ and $f(r)=r^{2}$ on $[0,R]$. The components $\S^{1}$ homogeneous of weight $m$ of the isoresonant potentials of the theorem have the following form :
\[V_{m}(r_{1}e^{i\alpha_{1}},\ldots,r_{k}e^{i\alpha_{k}},z)=\sum_{\overset{(\ell_{1},\ldots,\ell_{k})\in\Z^{k}}{\underset{i=1}{\overset{k}{\sum}}\ell_{i}p_{i}=-m}} W_{m,\ell_{1},\ldots,\ell_{k}}(\bar{x})e^{i\ell_{1}\alpha_{1}}\ldots e^{i\ell_{k}\alpha_{k}}\]
where $(r_{1}e^{i\alpha_{1}},\ldots,r_{k}e^{i\alpha_{k}},z)\in(\R^{2})^{k}\times\R^{n-2k}$, $\bar{x}\in\R^{n}/\S^{1}$ and the sum converges in infinite norm. 
\;\;

\item For the hyperbolic space $\H^{n}$, we can take the Poincar\'e model i.e. the unit ball of $\R^{n}$ centered at the origin with the metric $4(1-\parallel\! x \!\parallel^{2})^{-2}g_{euclid}$. The action of $\S^{1}$ on $\R^{n}$ described in the previous point induces an isometric action of $\S^{1}$ on $\H^{n}$. For the condition $C$, if $K$ is included in a ball of radius $R$ and centered at the origin (i.e. $0_{\R^{n}}$ with this model), then we take for $\widetilde{K}$ a ball of radius $\widetilde{R}>R$ and, in polar coordinates,  $\widetilde{g}=dr^{2}+f(r)d\omega^{2}$ with this time $f(r)=\sh^{2} r$ on $[0,R]$. The isoresonant potentials have the same form as in the Euclidean case. We recall (\cite{GZ1}) that if $n$ is odd $\Res(\Delta)=\emptyset$ and if $n$ is even $\Res(\Delta)=-\N$ and the multiplicity of the integer $-k$ is the multiplicity of $k(k+n-1)$ as eigenvalue of the Laplacian on the Euclidean sphere $\S^{n}$.

\;\;

\item Let us consider $\H^{n}$ with the model $\R^{+}_{\ast}\times\R^{n-1}$ with the corresponding coordinates $(x,y)$. We take for $X$ the hyperbolic cylinder $\H^{n}/\!<\!\gamma\!>$ where $\gamma$ is the isometry $w\rightarrow e^{\ell}w$. $\S^{1}$ acts on $X$ isometrically by $e^{i\theta}.[x,y]=[e^{\frac{\ell\theta}{2\pi}}x,
e^{\frac{\ell\theta}{2\pi}}y]$. We can see $X$ as $\R^{+}\times \S^{1}\times \S^{n-2}$ with the coordinates $(r,\theta,\omega)$, the metric $\d r^{2}+\ch^{2}r \d\theta^{2}+\sh^{2}r \d\omega^{2}$, and the $\S^{1}$ action is the trivial action on the factor $\S^{1}$. Every compact is included in a $K=[0,R]\times\S^{1}\times\S^{n-2}$, so we can take $\widetilde{K}=[0,\widetilde{R}]\times\S^{1}\times\S^{n-2}$ with the metric $\widetilde{g}=\d r^{2}+h_{1}(r)\d\theta^{2}+h_{2}(r)\d\omega^{2}$ where $h_{1}$ and $h_{2}$ are smooth on $[0,\widetilde{R}]$, constant near $\widetilde{R}$ and on $[0,R]$, $h_{1}(r)=\ch^{2}r$, $h_{2}(r)=\sh^{2}r$. The components $\S^{1}$ homogeneous of weight $m$ of the isoresonant potentials have the form :
\[V_{m}(r,\theta,\omega)=W_{m}(r,\omega)e^{-im\theta}.\] 
Here, according to \cite{GZ1}, we have $\Res(\Delta)=-\N+i\Z2\pi/\ell$ \cite{GZ1}.
\end{itemize}

\vskip 0.5cm

We recall that the isometric $\S^{1}$ action induces a decomposition of $L^{2}(X)$ according to isotypical subspaces :
\[L^{2}(X)=\bigoplus_{j\in\,\Z} L^{2}_{j}(X).\]

\noindent Let $P_{j} : L^{2}(X)\longrightarrow L^{2}_{j}(X)$ be the corresponding orthogonal projection.

The main idea of the proof is that if $V_{m}\in L^{\infty}(X)$ is $\S^{1}$ homogeneous of weight $m$ then it induces by multiplication a shift on these isotypical representations :
\[V_{m} : L^{2}_{j}(X)\longrightarrow L^{2}_{j+m}(X).\]

\subsection{Spectral lower bound for the Laplacian on homogeneous functions}

In the following we will need a spectral lower bound for the Laplacian on $\S^{1}$ homogeneous and compactly supported functions. As I have not read this result anywhere before, I also give an upper bound for the first eigenvalue. In \cite{AA}, it can be found a more precise discussion about this result.

\begin{proposition}
\label{prop:spectre}
Let $K$ be a compact manifold with boundary, with a $\S^{1}$ action and equipped with a metric $g$ such that $\S^{1}$ acts by isometries on $(K,g)$. We assume $g$ has a product form $\d\delta^{2}+h_{\partial K}$ in a neighbourhood of $\partial K$ with $\delta$ a $\S^{1}$ invariant function defining $\partial K$, and $h$ is independent of $\delta$. Then there are constants 
$C_{1}=C_{1}(K)>0$ and $C_{2}=C_{2}(K)>0$ such that for all $j\in\Z$, we have :
\[C_{1} j^{2}\leq\min \; \spec \; \Delta_{L^{2}_{j}(K)}\leq C_{2}<j>^{2}.\]
where $\Delta_{L^{2}_{j}(K)}$ is the Friedrichs selfadjoint extension in $L^{2}_{j}(K)$ of the Laplacian defined on $C^{\infty}_{c}(K)\cap L^{2}_{j}(K)$ and $<j>:=(1+\mid\!j\!\mid^{2})^{\frac{1}{2}}$.
\end{proposition} 

\begin{rem}
In the case where $K$ is a disk centered in $0$ in $\R^{2}$, we can apply this lemma including $K$ in a bigger disk $\widetilde{K}$ as explained before but we can also prove directly the lower bound because it is just an estimate of the first zero of Bessel functions.
\end{rem}

Proof : We begin with the lower bound. Consider two copies of $K$. We can identify their regular boundaries and get a compact closed manifold $M$. More precisely, there exists a collar neighbourhood $W$ of $\partial K$ diffeomorphic to $[0,\epsilon)\times\partial K$ by the diffeomorphism :
\begin{eqnarray*}
\psi : [0,\epsilon)\times\partial K &\rightarrow& W\\
(t,y)&\rightarrow& \psi_{t}(y),
\end{eqnarray*}
with $\psi_{t}$ the gradient flow of $\delta$ for the metric $g$. So, on the topological space $M=(K\sqcup K)/\partial K$, we can construct a differential atlas beginning with $\partial K\subset M$ which is included in a open set $[W]=(W\sqcup W)/\partial K$ diffeomorphic to $(-\epsilon,\epsilon)\times\partial K$ by
\begin{eqnarray*}
 (-\epsilon,\epsilon)\times\partial K &\rightarrow& [W]\\
(t,y)&\rightarrow& \left\{ \begin{array}{ll}
                   \psi_{t}(y),  &\textrm{if $t\geq 0$}\\
		   \psi_{-t}(y),  &\textrm{if $t\leq 0$}
                          \end{array} \right.
\end{eqnarray*}
On $[W]$ the $\S^{1}$ action is, via the previous diffeomorphism, the action on $\partial K$. The other charts are those in the interior of $K$.

As the metric $g$ has a product form in a neighbourhood of $\partial K$, it can be continue by symmetry on $\delta$. We get a smooth metric on $M$ and we still have an isometrical action of $\S^{1}$ on $M$. \\
Let $Y$ be the corresponding vector field which we will consider as a differential
operator of order $1$ ($Y.f(m)=-i\partial_{\theta}(f(e^{-i\theta}.m))_{\mid\theta=0}$). Another pseudo-differential operator of order $1$, on $M$, is $P:=\sqrt{\Delta_{M} +1}$ where $\Delta_{M}$ is the Laplacian on $(M,g)$. $P$ and $Y$ commute because $\S^{1}$ acts 
by isometries on $M$. We consider $Q:=P^{2}+Y^{2}$ whose principal symbol is $q(x,\xi)=\mid\!\xi\!\mid^{2}+(\xi(Y))^{2}$, $(x,\xi)\in T^{\ast}M$ ; 
so $Q$ is elliptic.

\;\;

Let $\Lambda\subset\R^{2}$ be the joint spectrum of $(P,Y)$, it is constituted
by the points $(\lambda_{k}^{P},\lambda_{k}^{Y})$ such that $P\phi_{k}=\lambda_{k}^{P}\phi_{k}$ and $Y\phi_{k}=\lambda_{k}^{Y}\phi_{k}$
where $(\phi_{k})$ is a orthonormal basis of $L^{2}(M)$. We note that the spectrum of $Y$ is equal to $\Z$ and we are looking for the minimum
of the first coordinates of points of $\Lambda$ whose second coordinate is $j$. We call this minimum $\lambda^{j}_{1}$.

\;\;

Let $p(x,\xi)=(\mid\!\xi\!\mid,\xi(Y))$ be the joint principal symbol of $P$ and $Y$. $p$ is an homogeneous function of degree $1$. Let 
$\Gamma$ be the linear cone $\Gamma=\R^{+}p(S(T^{\ast}M))$ where $S(T^{\ast}M)$ is the unit sphere bundle of $T^{\ast}M$. Moreover, we have 
$p(S(T^{\ast}M))=\{(1,\xi(Y)),\mid\!\xi\!\mid=1\}$. Now, as $P$ and $Y$ commute and $Q$ is elliptic, we can apply theorem $0.6$ in
\cite{CV} : if $C$ is a cone of $\R^{2}$ such that 
$C\cap\Gamma=\{0\}$ then $C\cap\Lambda$ is finite. Taking, for example, $C_{K}=\R\{(a,\mid\!Y\!\mid_{K}), \mid\! a\!\mid\leq\frac{1}{2}\}$ where $\mid\!Y\!\mid_{K}=\displaystyle\sup_{m\in K}\mid\!Y(m)\!\mid $, it means that there exists a constant $c:=\frac{1}{2\mid\!Y\!\mid_{K}}$ and $J\in\N$ such that, for all 
$j\in\Z$, $\mid\!j\!\mid\geq J$, we have $\lambda^{j}_{1}\geq c\mid\!j\!\mid$ and we can take a smaller $c$ and have $\lambda^{j}_{1}\geq c\mid\!j\!\mid$ for all $j\in\Z$. As $P=\sqrt{\Delta_{M}+1}$, there is another constant $c$, such that the minimum of the spectrum of the Laplacian on homogeneous functions of weight $j$ on $M$ is superior than $cj^{2}$.

\;\;

Moreover, the spectrum of $\Delta_{L^{2}_{j}(K)}$ with Dirichlet conditions is included in the spectrum
of $\Delta_{M}$ acting on $L^{2}_{j}(M)$. Indeed, let $I$ be the involution which exchanges the two copies of $K$ in
$M$, then the eigenfunctions of $\Delta_{M}$ which are odd for $I$ vanish on the image of $\partial K$ in $M$. So they correspond to
eigenfunctions of the Dirichlet Laplacian on $K$.\\

In order to prove the upper bound, we remark that
\[\min \; \spec \; \Delta_{L^{2}_{j}(K)}=\inf_{\phi\in L^{2}_{j}(K)\setminus\{0\}}\frac{\langle\Delta\phi ,\phi\rangle}{\parallel\!\phi\!\parallel^{2}}.\] 
So it's sufficient to construct, for $j$ large, one $\phi_{j}\in L^{2}_{j}(K)\setminus\{0\}$ such that $\langle\Delta\phi_{j} ,\phi_{j}\rangle\leq Cj^{2}\parallel\!\phi_{j}\!\parallel^{2}$ with $C$ independent of $j$. Let $\widehat{K}$ be the set of principal orbits of the action of $\S^{1}$ in $K$. The principal orbits are those for which the stability groups are the identity. $\widehat{K}$ is an open, connected and dense subset of $X$. We consider the principal $\S^{1}$ fibration $\widehat{K}\rightarrow \widehat{K}/\S^{1}$ and we take $U_{j}$, a $\S^{1}$ invariant open set of $\widehat{K}$, where this fibration is trivial. So $U_{j}$ is diffeomorphic to $(U_{j}/\S^{1})\times \S^{1}$ and we can take $U_{j}$ small enough to be sure that $U_{j}/\S^{1}$ is a coordinate patch and we denote the corresponding coordinates by $(\mathbf{y},\theta):=(y_{1},\ldots, y_{N}, \theta)$ where $N=\dim\ U_{j}/\S^{1}$. In those coordinates the metric has the form
\[g_{\mid U_{\phi_{j}}}=\sum_{k,\ell=1}^{N}a_{k,\ell}(\mathbf{y},\theta)\d y_{k}\d y_{\ell}+b(\mathbf{y}, \theta)\d\theta^{2}+\sum_{k=1}^{N}c_{k}(\mathbf{y}, \theta)\d y_{k}\d\theta,\]
with $a_{k,\ell}, b$ and $c_{k}$ smooth on $U_{j}$, and the Laplacian becomes
\begin{align*}
\Delta=& \sum_{k,\ell=1}^{N}A_{k,\ell}(\mathbf{y},\theta)\partial_{y_{k}}\partial_{y_{\ell}}+ B(\mathbf{y},\theta)\partial^{2}_{\theta}+\sum_{k=1}^{N}C_{k}(\mathbf{y},\theta)\partial_{\theta}\partial_{y_{k}}\\
&+\sum_{k=1}^{N}D_{k}(\mathbf{y},\theta)\partial_{y_{k}}+E(\mathbf{y},\theta)\partial_{\theta},
\end{align*}
where $A_{k,\ell}, B, C_{k}, D_{k}$ and $E$ are smooth on $U_{j}$.

\smallskip

We take $\phi_{j}(\mathbf{y},\theta)=\psi(\mathbf{y})e^{-ij\theta}$ with $\psi$ smooth and compactly supported in $U_{j}/\S^{1}$. So we have $\phi_{j}\in L^{2}_{j}(K)$ and 
\begin{align*}
\Delta\phi_{j}=& \big(\sum_{k,\ell=1}^{N}A_{k,\ell}(\mathbf{y},\theta)\partial_{y_{k}}\partial_{y_{\ell}}\psi-j^{2} B(\mathbf{y},\theta)\psi-ij\sum_{k=1}^{N}C_{k}(\mathbf{y},\theta)\partial_{y_{k}}\psi\\
&+\sum_{k=1}^{N}D_{k}(\mathbf{y},\theta)\partial_{y_{k}}\psi-ijE(\mathbf{y},\theta)\psi\big)e^{ij\theta},
\end{align*}
and
\begin{align*}
\langle\Delta\phi_{j},\phi_{j}\rangle=&-j^{2}\int_{\supp\,\phi_{j}}B(\mathbf{y},\theta)\mid\!\psi\!\mid^{2}\dvol(g)\\
&-ij\int_{\supp\,\phi_{j}}\sum_{k=1}^{N}C_{k}(\mathbf{y},\theta)(\partial_{y_{k}}\psi)\overline{\psi}+E(\mathbf{y},\theta)\mid\!\psi\!\mid^{2}\dvol(g)\\
&+\int_{\supp\,\phi_{j}}\sum_{k,\ell=1}^{N}A_{k,\ell}(\mathbf{y},\theta)(\partial_{y_{k}}\partial_{y_{\ell}}\psi)\overline{\psi}+\sum_{k=1}^{N}D_{k}(\mathbf{y},\theta)(\partial_{y_{k}}\psi)\overline{\psi} \,\, \dvol(g)
\end{align*}
so 
\[\langle\Delta\phi_{j},\phi_{j}\rangle=\mid\langle\Delta\phi_{j},\phi_{j}\rangle\mid\leq\Lambda_{1}(\psi)j^{2}\parallel\!\psi\!\parallel^{2}_{L^{2}(K,g)}+\Lambda_{2}(\psi)\mid\!j\!\mid+\Lambda_{3}(\psi),\]
where $\Lambda_{1}, \Lambda_{2}, \Lambda_{3}$ are positives constants which only depend of $\psi$. 

By fixing $\psi$ as $\parallel\!\phi_j\!\parallel=\parallel\!\psi\!\parallel$, we get a constant $C_{2}$ such that, for all $j\in\Z$,  
\[\langle\Delta\phi_{j},\phi_{j}\rangle\leq C_{2}<\!j\!>^{2}\parallel\!\phi_{j}\!\parallel^{2}. \ \ \square\]

\;\;\;

With the proposition \ref{prop:spectre} and the Condition $C$ we will prove the following :

\begin{lemme}
\label{lemme:majoration}
Let $\lambda\in D_{N}^{+}\setminus\Res(\Delta)$ and
$\chi\in C^{\infty}_{c}(X)$ be $\S^{1}$-invariant.
Then there is a constant $C=C(\lambda,\chi)>~0$ such that, for all $j\in\Z$, 
:
\[\parallel\chi \widetilde{R}_{0}(\lambda)P_{j}\chi\parallel\leq \frac{C}{1+j^{2}}\]
\end{lemme}

Proof : Since $\S^{1}$ acts by isometries on $X$, we have, for all $j$, $\Delta P_{j}=P_{j}\Delta$ and $P_{j}\widetilde{R}_{0}=\widetilde{R}_{0}P_{j}$. The fact that $\chi$ is $\S^{1}$ invariant also gives $\chi P_{j}=P_{j}\chi$.\\

We have 
\[\chi(\Delta-f(\lambda))\widetilde{R}_{0}(\lambda)\chi=\chi^{2},\]
so
\[(\Delta-f(\lambda))\chi \widetilde{R}_{0}(\lambda)P_{j}\chi=\chi^{2}P_{j}+[\Delta,\chi]\widetilde{R}_{0}(\lambda)\chi P_{j}.\]
then
\begin{equation}
\label{eq:const}
\parallel(\Delta-f(\lambda))\chi
\widetilde{R}_{0}(\lambda)P_{j}\chi\parallel \ \leq \ \parallel\chi^{2}P_{j}\parallel+\parallel [\Delta,\chi]\widetilde{R}_{0}(\lambda)\chi P_{j}\parallel \ \leq  C(\lambda,\chi).
\end{equation}

\;\;\;

Let $(\widetilde{K},\widetilde{g})$ the compact containing $K:=\supp\chi$ given by the Condition $C$. If $v\in L^{2}(X)$ then $u=\chi P_{j}\widetilde{R}_{0}(\lambda)\chi v$ is in $L^{2}_{j}(K,g)$ and, as $\supp\   u\subset K$ and $\widetilde{g}_{\mid\!K}=g_{\mid\!K}$, we also have $u\in L^{2}_{j}(\widetilde{K},\widetilde{g})$. In addition, $u$ is, at the same time, in the domain of the Dirichlet Laplacian on $\widetilde{K}$, of the Dirichlet Laplacian on $K$ and of the Laplacian on $X$, and we have
\[\Delta_{(\widetilde{K},\widetilde{g})}u=\Delta_{(K,g)}u=\Delta_{(X,g)}u.\]

Let $(\phi_{k})$ (depending on $\widetilde{K}$ and $j$) an orthonormal
basis of $L^{2}_{j}(\widetilde{K},\widetilde{g})$ constituted by eigenfunctions of $\Delta_{(\widetilde{K},\widetilde{g})}$. We denote $\mu_{k}(j,\widetilde{K})$ the
eigenvalue corresponding to $\phi_{k}$.
If we expand $u$ following this basis : $u=\underset{k}{\sum}u_{k}\phi_{k}$,
we have
\begin{equation*}
\big(\Delta_{(\widetilde{K},\widetilde{g})}-f(\lambda)\big)u = \sum_{k}\big(\mu_{k}(j,\widetilde{K})-f(\lambda)\big)u_{k}\phi_{k}
\end{equation*}
so that
\begin{equation*}
\parallel\big(\Delta-f(\lambda)\big)u\parallel^{2}=\sum_{k}\mid\mu_{k}(j,\widetilde{K})-f(\lambda)\mid^{2}\mid\!
u_{k}\!\mid^{2} \geq \sum_{k}\big(\mu_{k}(j,\widetilde{K})-\Re(f(\lambda))\big)^{2}\mid\!
u_{k}\!\mid^{2}.
\end{equation*}

Thus, using the proposition \ref{prop:spectre}, there exists a constant $C=C(\widetilde{K})>0$ such that
\[\forall k\in\N, \;  \forall j\in\Z  \;\; \mu_{k}(j,\widetilde{K})\geq C j^{2}.\]
We take $J$ in order to have $C J^{2}>\Re(f(\lambda))$, then for all $\mid\! j\!\mid\geq J$, we have
\begin{equation*}
\parallel\big(\Delta-f(\lambda)\big)u\parallel^{2} \geq \big(C j^{2}-\Re(f(\lambda))\big)^{2}\sum_{k}\mid\! u_{k}\!\mid^{2} = \big(C j^{2}-\Re(f(\lambda))\big)^{2} \parallel\! u\!\parallel^{2}.
\end{equation*}

Using this in the inequality (\ref{eq:const}) we get that for all $\mid\! j\!\mid\geq J$
\[\parallel\chi \widetilde{R}_{0}(\lambda)P_{j}\chi\parallel \leq \frac{C(\lambda,\chi)}{Cj^{2}-\Re(f(\lambda))}.\]
So there exists another constant $C>0$ such that , for all $\mid\! j\!\mid\geq J$, $\parallel\chi \widetilde{R}_{0}(\lambda)P_{j}\chi\parallel \leq Cj^{-2}$, and we can take a greater $C$ to have, for all $j\in\Z$,
\[\parallel\chi \widetilde{R}_{0}(\lambda)P_{j}\chi\parallel \leq\frac{C}{1+j^{2}}\ .\ \ \square\]

\subsection{Localization of resonances}

We begin the proof of the theorem \ref{thm : circular} by the inclusion $\Res(\Delta+V)\subset\Res(\Delta)$. First, we consider truncations in space of partial sums of $V$.

\subsubsection{Localization of resonances for the truncated partial sums of $V$}

We consider $S_{M}:=\underset{m=1}{\overset{M}{\sum}}V_{m}$ where $V_{m}$ is the component $\S^{1}$ homogeneous of weight $m$ of $V$. Let $\chi\in C^{\infty}_{c}(X)$ be invariant under the action of $\S^{1}$. In this part, our purpose is to show that $\Res(\Delta +\chi S_{M})\subset\Res(\Delta)$, on $D_{N}^{+}$.\\

For $\lambda\in D_{N}^{+}\setminus \Res(\Delta)$, we have
\[\big(\Delta+\chi S_{M}-f(\lambda)\big)\widetilde{R}_{0}(\lambda)\rho^{N} = \rho^{N}\big(I+\rho^{-N}\chi S_{M}\widetilde{R}_{0}(\lambda)\rho^{N}\big).\]
In addition, we have
\[\rho^{-N}\chi S_{M}\widetilde{R}_{0}(\lambda)\rho^{N}=\chi\rho^{-2N} S_{M}\rho^{N}\widetilde{R}_{0}(\lambda)\rho^{N}.\]
So thanks to the condition $A_{N,\rho}$, $\rho^{-N}\chi S_{M}\widetilde{R}_{0}(\lambda)\rho^{N}$ is an holomorphic family in $D_{N}^{+}\setminus\Res(\Delta)$, of compact operators such that
\[\parallel\!\rho^{-N}\chi S_{M}\widetilde{R}_{0}(\lambda)\rho^{N}\!\parallel<1\]
for $\mid\!\!\lambda\!\!\mid$ sufficiently large in $D_{N}^{+}$.

\noindent Then by the analytic Fredholm theory we get that $\big(I+\rho^{-N}\chi S_{M}\widetilde{R}_{0}(\lambda)\rho^{N}\big)^{-1}$ is meromorphic on $D_{N}^{+}\setminus \Res(\Delta)$ and we have the often called Lipmann-Schwinger equation which establishes the link between the resolvent of $\Delta+\chi S_{M}$ and these of the free Laplacian :
\begin{equation*}
\label{LS}
\rho^{N}\widetilde{R}_{\chi S_{M}}(\lambda)\rho^{N} = \rho^{N}\widetilde{R}_{0}(\lambda)\rho^{N}\big(I+\rho^{-N}\chi S_{M}\widetilde{R}_{0}(\lambda)\rho^{N}\big)^{-1}.
\tag{LS}
\end{equation*}
So if $\lambda_{0}$ is a pole of $\widetilde{R}_{\chi S_{M}}$ in $D_{N}^{+}\setminus\Res(\Delta)$, then $\lambda_{0}$ is a pole of 
$\big(I+\rho^{-N}\chi S_{M}\widetilde{R}_{0}(\lambda)\rho^{N}\big)^{-1}$ and still by Fredholm theory, there is a nontrivial $u\in L^{2}(X)$ such that
\[\big(I+\rho^{-N}\chi S_{M}\widetilde{R}_{0}(\lambda)\rho^{N}\big)u=0.\]

\;\;

\noindent We remark with the last equality that $\supp\ u\subset\supp\ \chi$. Let $\chi_{2}\in C^{\infty}_{c}(X)$ invariant under the action of $\S^{1}$ and such that $\chi_{2}=1$ on the support of $\chi$. If we denote $u_{j}:=P_{j}u\in L^{2}_{j}(X)$, we have
\[u_{j}= P_{j}\big(-\rho^{-N}\chi S_{M}\widetilde{R}_{0}(\lambda)\rho^{N}u\big) = P_{j}\big(-\rho^{-N}\chi S_{M}\chi_{2}\widetilde{R}_{0}(\lambda)\chi_{2}\rho^{N}u\big),\]
and by linearity
\[u_{j}=-\sum_{m=1}^{M}P_{j}\big(\rho^{-N}\chi V_{m}\chi_{2}\widetilde{R}_{0}(\lambda)\chi_{2}\rho^{N}u\big).\]

\;\;

However, each $V_{m}$ induces a shift on the isotypical representations : \[V_{m} : L^{2}_{j}(X)\rightarrow L^{2}_{j+m}(X),\] so we have
\begin{eqnarray*}
u_{j}&=&-\sum_{m=1}^{M}V_{m}P_{j-m}\big(\rho^{-N}\chi\chi_{2}\widetilde{R}_{0}(\lambda)\chi_{2}\rho^{N}u\big)\\
u_{j}&=&-\sum_{m=1}^{M}V_{m}\rho^{-N}\chi\chi_{2}\widetilde{R}_{0}(\lambda)\chi_{2}\rho^{N}P_{j-m}(u),
\end{eqnarray*}
where we have also used that the projections $P_{j-m}$ commute with $\widetilde{R}_{0}$, $\rho$, $\chi$ and $\chi_{2}$.

By hypothesis, for all $m$, $\parallel\!V_{m}\!\parallel_{\infty}\leq\underset{m'=1}{\overset{+\infty}{\sum}}\parallel\!V_{m'}\!\parallel_{\infty}<+\infty$, then, applying the lemma \ref{lemme:majoration} to $\chi_{2}\widetilde{R}_{0}\chi_{2}P_{j-m}$, we get a constant $C$ such that, for all $j\in\Z$,
\[\parallel u_{j}\parallel \leq \sum_{m=1}^{M}\frac{C}{1+(j-m)^{2}} \parallel u_{j-m}\parallel,\]
so, for all $j\in\Z$,
\[\parallel u_{j}\parallel \leq \epsilon_{j}\sum_{m=1}^{M}\parallel u_{j-m}\parallel,\]
where $\epsilon_{j}\rightarrow 0$ for $\mid\! j\!\mid\rightarrow +\infty.$

\;\;

\noindent Thus we can use the following lemma :

\begin{lemme}
\label{lemme:l2}
Let $(a_{j})_{j\in\Z}\in \ell^{1}(\Z)$ non-negative. If there is $M\in\N$ and, for all $j\in\Z$, $a_{j}\leq \epsilon_{j}\underset{m=1}{\overset{M}{\sum}}a_{j-m}$ with $\epsilon_{j}\rightarrow 0$ 
for $\mid\! j\!\mid\rightarrow +\infty$, then $a_{j}=0$ for all $j$.
\end{lemme}

\noindent Proof : Let $J'\leq 0$ such that, for all $j\leq J'$, $\epsilon_{j}\leq\frac{1}{M}$. Then, for all $j\leq J'$, we have $a_{j}\leq\frac{1}{M}\underset{m=1}{\overset{M}{\sum}}a_{j-m}$ and if we sum all these inequalities we get, denoting $S=\underset{j\leq J'}{\sum}a_{j}$,
\begin{align*}
S &\leq\frac{1}{M}\big((S-a_{J'})+(S-a_{J'}-a_{J'-1})+\ldots +(S-a_{J'}-\ldots -a_{J'-M+1})\big)\\
&=\frac{1}{M}\big(MS-Ma_{J'}-(M-1)a_{J'-1}-\ldots -a_{J'-M+1}\big),
\end{align*}
from which we deduce
\[0\leq -Ma_{J'}-(M-1)a_{J'-1}-\ldots -a_{J'-M+1},\]
and thus
\[a_{J'}=a_{J'-1}=\ldots=a_{J'-M+1}=0.\]
Moreover, as $\epsilon_{j}\rightarrow 0$ 
for $\mid\! j\!\mid\rightarrow +\infty$, there exists a constant $C$ such that, for all $j\in\Z$, $a_{j}\leq C\underset{m=1}{\overset{M}{\sum}}a_{j-m}$, so we have
\[\forall j\geq J'-M+1,   \,\,\,  a_{j}=0,\]
and we can make tighten $J'$ to $-\infty$ and finally we have $a_{j}=0$ for all $j\in\Z \ \ \ \square$

\;\;\;

We apply this lemma \ref{lemme:l2} to the sequence $\{\parallel\! u_{j}\!\parallel^{2}\}_{j}$. We get that $\parallel\!
u_{j}\!\parallel=0$ for all $j$ and thus $u\equiv 0$. This is in contradiction with the existence of a pole of $\widetilde{R}_{\chi S_{m}}$ in $D_{N}^{+}\setminus\Res(\Delta)$.

\;\;

\noindent Finally, for all $M$ and all $\chi\in C^{\infty}_{c}(X)$, invariant under the action of $\S^{1}$, $\Delta+\chi S_{M}$ has no resonance in $D_{N}^{+}\setminus\Res(\Delta)$, which can be expressed by 
\[\Res(\Delta+\chi S_{M})\subset\Res(\Delta), \ \ \textrm{in $D_{N}^{+}$}.\]

\subsubsection{Localization of resonances for the potential V}

Let us recall some results and notations about regularized determinant that will be needed in the following. (\cite{Y})

\begin{definition}
\label{def : schatten}
Let $\mathcal{H}$ be an Hilbert space. If $A:\mathcal{H}\longrightarrow\mathcal{H}$ is a compact operator, we define its
\textbf{singular values} $(s_{n}(A))_{n\in\N}$ as the eigenvalues of the selfadjoint operator $(A^{\ast}A)^{1/2}$. For $1\leq p<+\infty$, $\mathbf{\mathcal{S}_{p}}$, is the two-sided ideal of $\mathcal{L}(\mathcal{H})$ formed by operators $A$ for which the sum \[\parallel\! A\!\parallel_{p}^{p}=\sum_{n=0}^{\infty}s_{n}^{p}(A)\] is finite. 
\end{definition}

% \noindent Examples :
% 
% $\mathcal{S}_{1}=\{\textrm{trace class operators}\}$
% 
% $\mathcal{S}_{2}=\{\textrm{Hilbert-Schmidt operators}\}$
% 
% $\mathcal{S}_{\infty}:=\{\textrm{compact operators}\}$

\begin{definition}
For $A\in\mathcal{S}_{p}$ we define the \textbf{regularized determinant}, $\det_{p}$, by
\[\det_{p}(I+A)=\prod_{n=1}^{\infty}(1+\lambda_{n}(A))\exp\big(\sum_{k=1}^{p-1}\frac{(-1)^{k}}{k}\lambda_{n}^{k}(A)\big),\]
where the $(\lambda_{n}(A))_{n\in\N}$ are the eigenvalues of $A$.
\end{definition}

We give some properties of this determinant

\begin{proposition}
\label{prop:detreg}
\begin{itemize}

\item[1.] $A\longrightarrow\det_{p}(I+A)$ is continuous on
$(\mathcal{S}_{p},\parallel .\parallel_{p})$.

\;

\item[2.] If $z\longrightarrow A(z)$ is holomorphic in some domain of $\C$,
with values in $\mathcal{S}_{p}$, then $z\longrightarrow\det_{p}(I+A(z))$ is also holomorphic in the same domain.

\;

\item[3.] For $A\in\mathcal{S}_{p}$,  $I+A$ is invertible if and only if $\det_{p}(I+A)\neq 0$.
\end{itemize}
\end{proposition}

We have assumed that there exists $q$ such that for all $\lambda\in D_{N}^{+}\setminus\Res(\Delta)$, $\rho^{-(N+1)}V\widetilde{R}_{0}(\lambda)\rho^{N}$ is in a Schatten class $\mathcal{S}_{q}$, so $\rho^{-N}V\widetilde{R}_{0}(\lambda)\rho^{N}$ is in $\mathcal{S}_{q}$ too.

Let us first prove a preliminary fact. For all $\chi\in C^{\infty}_{c}(X)$ there exists $p$ such that $\chi\rho^{-N}\widetilde{R}_{0}(\lambda)\rho^{N}\in\mathcal{S}_{p}$ for all $\lambda\in D_{N}^{+}\setminus\Res(\Delta)$. To see this, take a compact $K$ with a smooth boundary and containing $\supp\chi$. Let $\Delta_{K}$ be the Dirichlet Laplacian on K, and $(\mu_{k})_{k\in\N}$ the eigenvalues of $(\Delta_{K}+1)^{-1}$. Then the Weyl's formula gives, when $k$ tends 
to $+\infty$,
\[\mu_{k}\sim \frac{(2\pi)^{2}}{(\omega_{n}\Vol(K))^{\frac{2}{n}}}k^{-\frac{2}{n}},\]
where $n=\dim X$ and $\omega_{n}$ is the volume of the unity ball in $\R^{n}$. Thus for $p>\frac{n}{2}$, $(\Delta_{K}+1)^{-1}\in\mathcal{S}_{p}$. Moreover, for all $\lambda\in D_{N}^{+}\setminus \Res(\Delta)$, $(\Delta_{K}+1)\chi\rho^{-N}\widetilde{R}_{0}(\lambda)\rho^{N}$ is a bounded operator in $L^{2}(X)$, and, as $\mathcal{S}_{p}$ is a two-sided ideal of $\mathcal{L}(L^{2}(X))$, we have
\[\chi\rho^{-N}\widetilde{R}_{0}(\lambda)\rho^{N}=(\Delta_{K}+1)^{-1}(\Delta_{K}+1)\chi\rho^{-N}\widetilde{R}_{0}(\lambda)\rho^{N}\in\mathcal{S}_{p}.\]

As $\mathcal{S}_{p_{1}}\subset\mathcal{S}_{p_{2}}$, for $p_{1}\leq p_{2}$, we can take the maximum of $p$ and $q$ and we still note it $q$, and get that $\chi\rho^{-N}\widetilde{R}_{0}(\lambda)\rho^{N}$ and $\rho^{-N}V\widetilde{R}_{0}(\lambda)\rho^{N}$ are both in $\mathcal{S}_{q}$.

The Lipmann-Schwinger equation, \eqref{LS}, with $V$ instead of $\chi S_{M}$ give
\[\rho^{N}\widetilde{R}_{V}(\lambda)\rho^{N}=\rho^{N}\widetilde{R}_{0}(\lambda)\rho^{N}\big(I+\rho^{-N}V\widetilde{R}_{0}(\lambda)\rho^{N}\big)^{-1}.\]
So thanks to the third point of the proposition \ref{prop:detreg} we have
\[\lambda\in\Res(\Delta +V)\cap D_{N}^{+}\setminus\Res(\Delta)\Longleftrightarrow \det_{q}\big(I+\rho^{-N}V\widetilde{R}_{0}(\lambda)\rho^{N}\big)=0.\]

\vskip 0.5cm

On $D_{N}^{+}\setminus \Res(\Delta)$, we define
\[F(V,\lambda):=\det_{q}\big(I+\rho^{-N}V\widetilde{R}_{0}(\lambda)\rho^{N}\big).\]
If there exists $\lambda_{0}\in \Res(\Delta +V)\setminus \Res(\Delta)$, then 
\[F(V,\lambda_{0})=0.\]
Let $\Gamma$ be a simple loop around $\lambda_{0}$ such that $\lambda_{0}$ is the only zero of $F(V,.)$ in the domain $U$ delimited by $\Gamma$, and
such that $\overline{U}\subset D_{N}^{+}\setminus \Res(\Delta)$. It is possible because, thanks to the second point of the proposition \ref{prop:detreg}, $F$ is holomorphic in $\lambda$ and so its zeros are isolated.\\

Let $\chi_{r}$ a smooth family of compactly supported and $\S^{1}$-invariant functions such that, $\parallel\!(\chi_{r}-1)\rho\!\parallel_{\infty}$ tends to $0$ when $r$ tends to $+\infty$. As we have assumed that $V\rho^{-(N+1)}\widetilde{R}_{0}(\lambda)\rho^{N}\in\mathcal{S}_{q}$ we can write, for all $\lambda\in\Gamma$,
\[\parallel\! \chi_{r}V\rho^{-N}\widetilde{R}_{0}(\lambda)\rho^{N}-V\rho^{-N}\widetilde{R}_{0}(\lambda)\rho^{N}\!\parallel_{q}\ \leq\ \parallel\!(\chi_{r}-1)\rho\!\parallel_{\infty}  \;  \parallel\!V\rho^{-(N+1)}\widetilde{R}_{0}(\lambda)\rho^{N}\!\parallel_{q}.\]
So, when $r$ tends to $+\infty$, $\chi_{r}V\rho^{-N}\widetilde{R}_{0}(\lambda)\rho^{N}$ tends to $V\rho^{-N}\widetilde{R}_{0}(\lambda)\rho^{N}$ in $\mathcal{S}_{q}$ uniformly on $\Gamma$. So, with the first point of the proposition \ref{prop:detreg}, $F(\chi_{r}V,\lambda)\rightarrow F(V,\lambda)$ uniformly on $\Gamma$.
From that, there exists $r_{0}$ such that for all $r>r_{0}$ and for all $\lambda\in\Gamma$ we have
\[\mid\! F(\chi_{r}V,\lambda)-F(V,\lambda)\!\mid  <  \mid \!F(V,\lambda)\!\mid.\]
So, by Rouché's theorem, $F(\chi_{r}V, .)$ has the same number of zeros, in $U$, as $F(V, .)$.\\

In the same way, fixing $r>r_{0}$, using $\chi_{r}\rho^{-N}\widetilde{R}_{0}(\lambda)\rho^{N}\in\mathcal{S}_{q}$ we can write
\[\parallel\! \chi_{r}S_{M}\rho^{-N}\widetilde{R}_{0}(\lambda)\rho^{N}-\chi_{r}V\rho^{-N}\widetilde{R}_{0}(\lambda)\rho^{N}\!\parallel_{q}\;\leq\;\parallel\!S_{M}-V\!\parallel_{\infty}  \;  \parallel\!\chi_{r}\rho^{-N}\widetilde{R}_{0}(\lambda)\rho^{N}\!\parallel_{q}.\]
So using the fact that by hypothesis, $\parallel\! S_{M}-V\!\parallel_{\infty}$ tends to $0$ when $M$ tends to $\infty$, we have $F(\chi_{r}S_{M},\lambda)\rightarrow F(\chi_{r}V,\lambda)$ uniformly on $\Gamma$ and we can use the Rouch\'e's theorem again.

In conclusion, there exist $r$ and $M$ such that $F(\chi_{r}S_{M}, .)$ has the same number of zeros, in $U$, as $F(V, .)$. It means that $\Delta+\chi_{r}S_{M}$ has a resonance in the domain $U\subset D_{N}^{+}\setminus\Res(\Delta)$ which is in contradiction with the previous part. In conclusion, on $D_{N}^{+}$,
\[\Res(\Delta+V)\subset\Res(\Delta).\]

\begin{rem}
In \cite{C3}, Christiansen proves the inclusion, $\Res(\Delta+V)\subset\Res(\Delta)$, without using the shift created by the potential $V$ on the isotypical components of $L^{2}(X)$ but with regularized determinant and an hypothesis of analycity : $W(z):=\underset{m=1}{\overset{\infty}{\sum}}z^{m}V_{m}$ should be holomorphic in a domain of $\C$ containing the closed disc of center $0$ and radius $1$.
\end{rem}

\;\;

\subsection{Persistence of resonances}

In order to achieve the proof of theorem \ref{thm : circular}, we have to show that the points in $\Res(\Delta)\cap D_{N}^{+}$ are also resonances of $\Delta+V$
with the same multiplicity. To make this, we will use the Agmon's perturbation theory of resonances \cite{A}.

\;\;

Let $\lambda_{0}\in D_{N}^{+}$ be a resonance of $\Delta$ with multiplicity $m$. Let $U\subset D_{N}^{+}$ with smooth boundary $\Gamma$ such that $\overline{U}\cap\Res(\Delta)=\{\lambda_{0}\}$. If $V$ satisfies the hypothesis of the theorem \ref{thm : circular} then, for all $t\geq 0$, $tV$ satisfies these hypothesis too. So we can apply the result of the previous part : for all $t\geq 0$, $\Res(\Delta+tV)\subset\Res(\Delta)$ and thus
\[\Res(\Delta+tV)\cap U\subset\{\lambda_{0}\}.\]
Let $E:=\{t_{0}\geq 0\ ;\ \forall t\in[0,t_{0}],\ \Res(\Delta+tV)\cap U=\{\lambda_{0}\} \ \ \textrm{with the same multiplicity}\  m\}$ ; we are going to prove by connexity that it is in fact equal to $[0,+\infty[$. First it is not empty because by definition of $\lambda_{0}$, $0\in E$.\\

We take $t_{0}\in E$, we want to prove that there exists $\delta>0$ such that $]t_{0}-\delta, t_{0}+\delta[\subset E$. Following the theory of Agmon (\cite{A}), we begin with the definition of the Banach space
\[B_{\Gamma}=\{f\in \rho^{-N}L^{2}(X) \; ; \; f=g+\int_{\Gamma}\widetilde{R}_{t_{0}V}(\xi)\Phi(\xi)\d\xi,     \;\;  g\in \rho^{N}L^{2}(X), \Phi\in C(\Gamma,\rho^{N}L^{2}(X))\},\]
where $C(\Gamma,\rho^{N}L^{2}(X))$ is the space of continuous functions on $\Gamma$ with values in $\rho^{N}L^{2}(X)$. On the space $B_{\Gamma}$ we take the norm
\[\parallel\!f\!\parallel_{B_{\Gamma}}=\inf_{g,\Phi}(\parallel\!g\!\parallel_{\rho^{N}L^{2}(X)} +\parallel\!\Phi\!\parallel_{C(\Gamma,\rho^{N}L^{2}(X))}),\]
where the infimum is taken among all the $g\in \rho^{N}L^{2}(X)$ and the $\Phi\in C(\Gamma,\rho^{N}L^{2}(X))$ such that $f=g+\int_{\Gamma}\widetilde{R}_{t_{0}V}(\xi)\Phi(\xi)\d\xi$. On this Banach space, we can define, still following Agmon, the operator $(\Delta+t_{0}V)^{\Gamma} : \mathcal{D}\big((\Delta+t_{0}V)^{\Gamma}\big)\rightarrow B_{\Gamma}$. It's a restriction of $\Delta+t_{0}V$ in the following sens : 
\[(\Delta+t_{0}V)^{\Gamma}u=(\overline{\Delta+t_{0}V})u, \ \  u\in \mathcal{D}\big((\Delta+t_{0}V)^{\Gamma}\big), \]
where $\overline{\Delta+t_{0}V}$ is the closure of the operator $\Delta+t_{0}V$ view as an operator densely defined in $\rho^{-N}L^{2}(X)$.\\
Agmon proves that $(\Delta+t_{0}V)^{\Gamma}$ has a discrete spectrum in $U$ which is exactly the set of the poles of $\widetilde{R}_{t_{0}V}$, i.e. the resonances of $\Delta+t_{0}V$, with the same multiplicities.\\

Next, with the condition $B_{N,\rho}$, the family $tV$ verifies all the hypothesis in order to apply the part "perturbation" of the paper \cite{A}. We perturb $\Delta+t_{0}V$ by $tV$. So there exists $\delta>0$ such that, for all $t\in]-\delta,\delta[$, we can define in $B_{\Gamma}$ the operator $(\Delta+t_{0}V+tV)^{\Gamma}$. Moreover, for all $t\in]-\delta,\delta[$, $(\Delta+t_{0}V+tV)^{\Gamma}$ has a discrete spectrum in $U$ which is exactly the set of the poles of $\widetilde{R}_{t_{0}V+tV}$ with the same multiplicities.\\

Now, our problem becomes a problem of eigenvalues. Using the Kato's perturbation theory of eigenvalues (\cite{K}), we know that, maybe taking a smaller $\delta$, the eigenvalues of $(\Delta+t_{0}V+tV)^{\Gamma}$ in $U$ are continuous for all $t\in]-\delta,\delta[$. As these eigenvalues are also resonances of $\Delta+t_{0}V+tV$, $\lambda_{0}$ is the unique possibility. So $\lambda_{0}$ is the unique eigenvalue in $U$ of $(\Delta+t_{0}V+tV)^{\Gamma}$ for all $t\in]-\delta,\delta[$ with constant multiplicity. Therefore, thanks to the parallel established before, $\lambda_{0}$ is the unique resonance in $U$ of $\Delta+t_{0}V+tV$ for all $t\in]-\delta,\delta[$ with constant multiplicity. It signifies $]t_{0}-\delta, t_{0}+\delta[\subset E$ and so $E$ is an open set.

We can prove that $E$ is also a closed set doing the same proof with the complementary set of $E$. If $t_{0}$ is not in $E$, then $\lambda_{0}$ is a resonance of $\Delta+t_{0}V$ with a multiplicity not equal to $m$ (it can be $0$). Perturbing this operator by $tV$ and using the Agmon's correspondence, we can prove that $\lambda_{0}$ is a resonance of $\Delta+tV$ with a multiplicity not equal to $m$ for all $t$ in a neighbourhood of $t_{0}$.\\

In conclusion, $E=[0,+\infty[$ and we can take $t_{0}=1$ to obtain, in $U$, $\Res(\Delta+V)=\Res(\Delta)$ with the same multiplicity. To finish, we have to do the same work in the neighbourhood of any resonance of the free Laplacian. This completes the proof of the theorem \ref{thm : circular}.

\;\;\;

\subsection{An example where the order of resonances grows}

The isoresonant potentials introduced in the theorem \ref{thm : circular} can't be detected only observing the set of resonances and their multiplicities. We can wonder if their existence can be seen through the order of the resonances. We are going to prove, in an example, that there exist potentials verifying the theorem \ref{thm : circular} which change the order of resonances.\\

We consider the hyperbolic plane $\H^{2}$ with the model $\R^{+}\times\S^{1}$, the coordinates $(r,\theta)$ and the metric $g=\d r^{2}+\sh (r)^{2}\d\theta^{2}$. We have already said that the resonances of the free Laplacian are all the negative integers and the multiplicity of $-k$, $k\in\N$, is $2k+1$ (see \cite{GZ1}). Moreover the order of all these resonances is $1$. We denote $\mathcal{F}:=\{V_{m}(r)e^{im\theta} ; m\in\Z\setminus\{0\}, V_{m}\in L^{\infty}_{c}(\R^{+})\}$ which is a family of isoresonant potentials by the theorem \ref{thm : circular}. 

\begin{proposition}
\label{prop:ordre}
Let $\H^{2}$ be the hyperbolic plane and $k$ a non negative integer. There exists a potential $V\in\mathcal{F}:=\{V_{m}(r)e^{im\theta} ; m\in\Z\setminus\{0\}, V_{m}\in L^{\infty}_{c}(\R^{+})\}$ such that $-k$ is a resonance of $\Delta+V$ with an order strictly greater than $1$.
\end{proposition}
Proof : let $k\in\N\setminus\{0\}$, we suppose, ad absurdum, for all $V\in\mathcal{F}$, $-k$ is a resonance of order $1$ of $\Delta+V$. \\
\indent For all $V\in\mathcal{F}$, the resolvent $(\Delta+V-\lambda(1-\lambda))^{-1}$ has a meromorphic continuation $\widetilde{R}_{V}$ on $D^{+}_{N}=\{\lambda\in\C\ ;\ \Re\lambda>\frac{1}{2}-N\}$ as an operator from $B_{0}:=\rho^{N}L^{2}(\H^{2})$ to $B_{1}:=\rho^{-N}L^{2}(\H^{2})$ where $\rho$ is a boundary defining function of a compactification of $\H^{2}$. We take $N$ sufficiently large to have $-k\in D_{N}^{+}$. $B_{0}$ and $B_{1}$ are dual thanks to the non degenerate symmetric form :
\[\langle u,v\rangle=\int_{\H^{n}}u v \,\dvol(g).\]
We remark, for all $t\in\R$ and all $V\in\mathcal{F}$, we have $tV\in\mathcal{F}$. With our hypothesis, for $\lambda$ in a neighbourhood of $-k$, we have
\[\widetilde{R}_{tV}(\lambda)=(\lambda+k)^{-1}S(t)+H(t,\lambda),\] 
where $H(t,.)$ is holomorphic with values in $\mathcal{L}(B_{0},B_{1})$ and $S(t)\in\mathcal{L}(B_{0},B_{1})$ has a finite rank.

We apply the Agmon's perturbation theory of resonances. Consider a domain $U\subset D^{+}_{N}$ with smooth boundary $\Gamma$ such that $\overline{U}\cap\Res(\Delta)=\{-k\}$. We have the corresponding Banach space,
\[B_{\Gamma}=\{f\in B_{1} \; ; \; f=g+\int_{\Gamma}\widetilde{R}_{0}(\xi)\Phi(\xi)\d\xi,     \;\;  g\in B_{0}, \Phi\in C(\Gamma,B_{0})\},\]
with $B_{0}\subset B_{\Gamma}\subset B_{1}$.\\
\indent Then there exists $\delta>0$ such that, for all $V\in\mathcal{F}$ and all $t\in]-\delta,\delta[$, we can define the operators $(\Delta+tV)^{\Gamma}$ in $B_{\Gamma}$ and their resolvents $R_{tV}^{\Gamma}$. Thanks to \cite{A}, we know that $(\Delta+tV)^{\Gamma}$ has a discrete spectrum in $U$ which correspond to the  resonances of $\Delta+tV$ in $U$ with the same multiplicities and orders. So for $\lambda$ near $-k$ in $U$ we have
\begin{equation}
\label{ordre1}
R_{tV}^{\Gamma}(\lambda)=(\lambda+k)^{-1}S^{\Gamma}(t)+H^{\Gamma}(t,\lambda),
\end{equation}
with $H^{\Gamma}(t,.)$ holomorphic with values in $\mathcal{L}(B_{\Gamma})$ and $S^{\Gamma}(t)\in\mathcal{L}(B_{\Gamma})$ of finite rank.
Still following \cite{A} we know that $S(t)$ and $S^{\Gamma}(t)$ have the same range and they coincide on $B_{0}$.

Let $V\in\mathcal{F}$, for all $t\in]-\delta,\delta[$ and $\phi\in B_{\Gamma}$ we define $\psi(t):=S^{\Gamma}(t)\phi$. From (\ref{ordre1}) we obtain for all $t\in]-\delta,\delta[$,
\[((\Delta+tV)^{\Gamma}+k(k+1))\psi(t)=0.\]
$\psi(t)$ is derivable in $t$ like $S^{\Gamma}(t)$ (because $S^{\Gamma}(t)=\frac{1}{2\pi i}\int_{\Gamma}(\Delta^{\Gamma}+tV-\lambda(1-\lambda))^{-1}d\lambda$), so we can derivate the last equality at $t=0$ and get
\[V\psi(0)+(\Delta^{\Gamma}+k(k+1))\psi'(0)=0.\]
Compose this new equality with $S^{\Gamma}(0)$, using, $S^{\Gamma}(0)(\Delta^{\Gamma}+k(k+1))=(\Delta^{\Gamma}+k(k+1))S^{\Gamma}(0)=0$, because 
\[S^{\Gamma}(0)=\frac{1}{2\pi i}\int_{\Gamma}(\Delta^{\Gamma}-\lambda(1-\lambda))^{-1}d\lambda,\]
and the fact that $-k(k+1)$ is an eigenvalue of order $1$ of $\Delta^{\Gamma}$, we obtain 
\[S^{\Gamma}(0)V\psi(0)=0.\]
As $\psi(0)\in\Ran S^{\Gamma}(0)=\Ran S(0)$, there exists $f_{0}\in B_{0}$ such that $\psi(0)=S(0)f_{0}=S^{\Gamma}(0)f_{0}$. Moreover
$S^{\Gamma}(0)V\psi(0)\in B_{\Gamma}\subset B_{1}$, so we can evaluate :
\[\langle S^{\Gamma}(0)V\psi(0) , f_{0}\rangle=0.\]
\indent Next, as the Laplacian is a real operator, it is symmetric for $\langle.,.\rangle$ thus the resolvent, $\widetilde{R}_{0}(\lambda)$, is symmetric too, first for $\Re(\lambda)>\frac{1}{2}$ and after in all $D^{+}_{N}$ by analytic continuation. In conclusion
$S^{\Gamma}(0)$ is symmetric for $\langle.,.\rangle$. So
\[\langle V\psi(0) , S^{\Gamma}(0)f_{0}\rangle=\langle V\psi(0) , \psi(0)\rangle=0,\]
and we have that equality for all $V\in\mathcal{F}$.

Thus, for all $m\in\Z\setminus\{0\}$ and all $V_{m}\in L^{\infty}_{c}(\R^{+})$ we have
\[\int_{0}^{2\pi}e^{im\theta}\int_{\R^{+}}V_{m}(r)\psi(0)^{2}(r,\theta) \dvol(g)=0.\] 
This implies that, for all the resonant states $\psi(0)$ of the free Laplacian, $\psi(0)^{2}$ does not depend on $\theta$. But, considering the expression of the hyperbolic Laplacian and taking its decomposition corresponding to  $\underset{\ell\in\Z}{\bigoplus}(L^{2}(\R^{+}, \sh r \d r)\otimes e^{i\ell\theta})$, we have resonant states of the form $\psi(0)(r,\theta)=\psi_{\ell}(r)e^{i\ell\theta}$ where $\mid\!\ell\!\mid\leq k$ and $\psi_{\ell}$ are hypergeometric functions (see the annexe of \cite{GZ2}). Then for $\ell\neq 0$, $\psi(0)^{2}$ depend on $\theta$ : we have our contradiction.

Finally there exists $V\in\mathcal{F}$ such that $-k$ is a resonance of $\Delta+V$ of order strictly greater than $1$.\ \ $\square$

\section{$SO(n)$ symmetries}

This time we consider an isometric action of $SO(n)$ on a complete Riemannian manifold $(X,g)$ of dimension $n\geq3$. Contrary to the case $\S^{1}$, $SO(n)$ is not commutative, so we don't have a simple description of the isotypical components. To have a shift we add an hypothesis~ :

\vskip 0.3cm

\textbf{\underline{Condition $D$}} : The isometric action of $SO(n)$ on $(X,g)$ has a fixed point $O$ and the polar coordinates with pole $O$ define a diffeomorphism from $X\setminus\{O\}$ to $\R^{+}\setminus\{0\}\times\S^{n-1}$. 

\vskip 0.3cm

With this condition, in the polar coordinates the metric $g$ becomes :
\[dr^{2}+f(r)\d\omega^{2}, \;\; (r,\omega)\in\R^{+}\times\S^{n-1}, \]
where $\d\omega^{2}$ is the metric on the $(n-1)$-sphere in $\R^{n}$. For example with $f(r)=r^{2}$ we have the Euclidean space and, with
$f(r)=\sh(r)^{2}$, the hyperbolic space. If $f$ is independent of $r$ outside a compact then $(X,g)$ is a manifold with a cylindrical end of section $\S^{n-1}$.

With the condition $D$ we have
\[L^{2}(X)=\bigoplus_{k\in\N}L^{2}(\R^{+})\otimes H^{k},\] 
where $H^{k}=\Ker(\Delta_{\S^{n-1}}-k(k+n-2))$, $k\in\N$, be the eigenspaces of the Laplacian on $\S^{n-1}$. The action of $SO(n)$ on $X$ induces a representation of $SO(n)$ on $L^{2}(X)\simeq L^{2}(\R^{+})\otimes L^{2}(\S^{n-1})$ which only acts on the factor $L^{2}(\S^{n-1})$, so on the $H^{k}$. Moreover the restriction of this representation to each $H^{k}$ is irreducible (cf \cite{BGM}). The shift that we will use in order to construct isoresonant potentials, will appear on these $H^{k}$.

\subsection{Representation}

The group action of $SO(n)$ on $L^{2}(X)$ induces an action of its Lie algebra,  $\mathfrak{so}_{n}$. We can describe this action with the following operators, 
\[D_{\xi}f(x):=\frac{\d}{\d t}f(e^{-t\xi}.x)_{\mid t=0}, \;\;\; \xi\in\mathfrak{so}_{n}, \ f\in L^{2}(X), \; x\in X.\]
We consider the complexification of the Lie algebra $\mathfrak{so}_{n}$, $\mathfrak{g}:=\mathfrak{so}_{n}^{\C}=\mathfrak{so}_{n}+i\mathfrak{so}_{n}$. We choose $\mathfrak{h}$ one of the Cartan subalgebras of $\mathfrak{g}$ i.e. one of the maximal Abelian subalgebras of $\mathfrak{g}$. Let us describe $\mathfrak{h}$ as a subalgabra of $\mathfrak{gl}(\C^{n})$. $\mathfrak{h}$ is the Lie algebra whose basis is $(\zeta_{k})_{1\leq k\leq p}$ where $p$ is the integer part of $\frac{n}{2}$ and $\zeta_{k}$ has all its entries null except the $k^{th}$ $2\times 2$-block which is 
\[\begin{pmatrix}
0 & i \cr
-i & 0 \cr
\end{pmatrix}.\]
Take $(\omega_{k})$ the dual basis of $(\zeta_{k})$ in $\mathfrak{h}^{\ast}$.

Remember that $\mathfrak{g}$ acts on itself by the adjoint representation :
\[\ad(Y) : Z\rightarrow [Y,Z], \;\;\; (Y,Z)\in\mathfrak{g}^{2}.\]
We consider the following scalar product, 
\[\langle Y,Z\rangle=\Tr(\ad(\overline{Y})\circ\ad(Z)), \ \ (Y,Z)\in\mathfrak{g}^{2}\]
where the conjugation is defined by $\overline{U+iV}=-U+iV$ with real $U$ and $V$. So, for all $Y,Z\in\mathfrak{g}$, $\overline{[Y,Z]}=-[\overline{Y},\overline{Z}]$. With this we remark that, for all $\xi\in\mathfrak{h}$, we have $\overline{\xi}=\xi$ and thus $\ad(\xi)$ is selfadjoint (cf \cite[p.177]{Si}). So $\{\ad(\xi)\; ;\; \xi\in\mathfrak{h}\}$ is a family of selfadjoint operators on $\mathfrak{g}$ which commute together. We can simultaneously diagonalize them and decompose $\mathfrak{g}$ according to the eigenspaces.

We obtain $\mathfrak{g}=\mathfrak{h}\oplus\bigoplus\mathfrak{g}_{\alpha}$ where the sum is over a finite set of $\alpha\in\mathfrak{h}^{\ast}$ which are the \textit{roots} of $\mathfrak{g}$ and we denote $\mathfrak{g}_{\alpha}:=\{X\in\mathfrak{g}\ ;\ \ad(\xi)(X)=\alpha(\xi)X ,  \;\; \forall\xi\in\mathfrak{h}\}$ which are the \textit{root spaces} (they are all one dimensional cf \cite[p.180]{Si}). Let $\Lambda\subset\mathfrak{h}^{\ast}$ the integer lattice generated by the roots. In $\Lambda$ we choose a lexicographical order ``$\succeq$'', choosing $\omega_{1}\succeq\ldots\succeq\omega_{p}$. Then we denote $\mathfrak{g}_{+}:=\underset{\alpha\succ 0}{\bigoplus}\mathfrak{g}_{\alpha}$ (respectively $\mathfrak{g}_{-}:=\underset{\alpha\prec 0}{\bigoplus}\mathfrak{g}_{\alpha}$) the subalgebra of $\mathfrak{g}$ generated by root spaces with positive root (respectively negative). So we have $\mathfrak{g}=\mathfrak{h}\oplus\mathfrak{g}_{+}\oplus\mathfrak{g}_{-}$. For a general theory  see \cite[chapter VIII]{Si}.

We come back to the irreducible representations $H^{k}$. It have a decomposition according to the action of $\mathfrak{h}$ :
\[H^{k}=\bigoplus_{\omega_{min}^{k}\preceq\omega\preceq\omega_{max}^{k}}H^{k}_{\omega},\]
where the sum is over a finite set of $\mathfrak{h}^{\ast}$, and these $\omega$ are the \textit{weights} of $H^{k}$ and the corresponding \textit{weight spaces}, $H^{k}_{\omega}$, are defined by $H^{k}_{\omega}=\{f\in H^{k} ;
D_{\xi}f=\omega(\xi)f, \,\, \forall\xi\in\mathfrak{h}\}$. 

We will need the following lemma,

\begin{lemme}
\label{liedecal}
If $f\in H^{k}_{\omega}$ and if $\xi\in\mathfrak{g}_{\alpha}$, then $D_{\xi}f\in H^{k}_{\omega+\alpha}$.
\end{lemme}
Proof : for all $\zeta\in\mathfrak{h}$, we have 
\[D_{\zeta}(D_{\xi}f)=D_{\xi}(D_{\zeta}f)+D_{[\zeta,\xi]}f.\]
But $[\zeta,\xi]=\ad(\zeta)(\xi)=\alpha(\zeta)\xi$ because $\xi\in\mathfrak{g}_{\alpha}$, and $D_{\zeta}f=\omega(\zeta)f$ by definition of $H^{k}_{\omega}$. So
\[D_{\zeta}(D_{\xi}f)=\omega(\zeta)D_{\xi}f+\alpha(\zeta)D_{\xi}f=(\omega+\alpha)(\zeta)(D_{\xi}f). \;\; \square\]

We define particular vectors in the $H^{k}$, $k\in\N$, which will be used to create the necessary shift.

\begin{definition}
A nonzero vector $v\in H^{k}$ is a \textit{highest weight vector} if it is an eigenvector for the action of all the $D_{\xi}$, $\xi\in\mathfrak{h}$, and if it is in the kernel of all the $D_{\xi}$, $\xi\in\mathfrak{g}^{+}$.
\end{definition}

As $\mathfrak{g}$ is semi-simple and $H^{k}$ is an irreducible representation of it, there is an unique highest weight vector up to scalar : we note it $v_{max}^{k}$. In fact $H^{k}_{\alpha_{max}^{k}}$ is one dimensional,  generated by $v_{max}^{k}$. With our choice $v_{max}^{k}$ can be calculated explicitly (see \cite{AA}) :
\[v_{max}^{k}\circ\phi\ (x_{1},\ldots,x_{n})=(x_{1}+ix_{2})^{k},\]
where $\phi :\R^{+}\setminus\{0\}\times\S^{n-1}\rightarrow X\setminus\{O\}$ is the diffeomorphism of condition $D$ and $(x_{1},x_{2},\ldots,x_{n})$ are the standard coordinates of $\R^{n}$ restricted to $\S^{n-1}$.

\vskip 0.3cm

We will need the following lemma,

\begin{lemme}
\label{lemmeGu}
\[H^{k}_{\omega_{max}^{k}}=H^{k}\cap\big(\bigcap_{\xi\in\mathfrak{g}_{+}} \Ker D_{\xi}\big).\]
\end{lemme}

Proof : the first inclusion $H^{k}_{\omega_{max}^{k}}\subset H^{k}\cap\big(\underset{\xi\in\mathfrak{g}_{+}}{\bigcap} \Ker D_{\xi}\big)$ is the definition of a highest weight vector. 

Let $u\in H^{k}\cap\big(\underset{\xi\in\mathfrak{g}_{+}}{\bigcap} \Ker D_{\xi}\big)$, so $u\in H^{k}=\!\underset{\omega_{min}^{k}\preceq\omega\preceq\omega_{max}^{k}}{\bigoplus}\!H^{k}_{\omega}$ and we write $u=\!\underset{\omega_{min}^{k}\preceq\omega\preceq\omega_{max}^{k}}{\sum}\!u_{\omega}$ with $u_{\omega}\in H^{k}_{\omega}$. For all $\xi\in\mathfrak{g}_{\beta}$ with $\beta\succ 0$, we have, thanks to the lemma \ref{liedecal}, $D_{\xi}u_{\omega}\in H^{k}_{\alpha+\beta}$. By hypothesis, we have 
\[D_{\xi}u=\!\sum_{\omega_{min}^{k}\preceq\omega\preceq\omega_{max}^{k}}\!D_{\xi}u_{\omega}=0,\]
and, as the previous sum is direct, we get for all $\omega$ :
\[\forall \xi\in\mathfrak{g}_{+}  \;\; D_{\xi}u_{\omega}=0.\] 
Moreover, for all $\omega$, by the definition of $H^{k}_{\omega}$, $u_{\omega}$ is an eigenvector for the action of all the $D_{\xi}$, $\xi\in\mathfrak{h}$. Thus by the definition of a highest weight vector, we have $u_{\omega}=0$ except for $u_{\omega_{max}^{k}}$ and in conclusion $u\in H^{k}_{\omega_{max}^{k}} \ \ \square$

\subsection{Isoresonant potentials} 

 Let 
\[L^{2}(X)^{+}=\bigoplus_{k\in\N}L^{2}(\R^{+})\otimes H^{k}_{\omega_{max}^{k}}\]
and note $P_{k}$ the corresponding projections into $L^{2}(\R^{+})\otimes H^{k}_{\omega_{max}^{k}}$.

\vskip 0.3cm

\begin{thm}
\label{thmSOn}
Let $(X,g)$ a Riemannian manifold of dimension $n\geq 3$, with an isometric action of $SO(n)$, verifying condition $D$. We assume that we have condition $A_{N,\rho}$ for some $N>0$ with a function $\rho$ invariant under the action of $SO(n)$.

Let $V$ be the potential,
\[V=\sum_{k=1}^{\infty}V_{k}\] 
where $V_{k}\in L^{2}(\R^{+})\otimes H^{k}_{\omega_{max}^{k}}$ and $\underset{k}{\sup}\parallel\!V_{k}\!\parallel_{\infty}<+\infty$. If $V$ verifies condition $B_{N,\rho}$, and for all $\lambda\in D_{N}^{+}\setminus\Res(\Delta)$, $\rho^{-(N+1)}V\widetilde{R}_{0}(\lambda)\rho^{N}$ is in a Schatten class $\mathcal{S}_{q}$, $q\in\N\setminus\{0\}$,\\

\noindent then, in $D_{N}^{+}$, $\Res(\Delta +V)=\Res(\Delta)$ with the same multiplicities.

\end{thm}

The Euclidean space $\R^{n}$, the hyperbolic space $\H^{n}$, asymptotically hyperbolic spaces and manifolds with asymptotically cylindrical ends with an action of $SO(n)$, are examples where this theorem can be applied.

\begin{rem}
If $X$ has an isometric action of $SO(n)$ it has also an isometric action of $\S^{1}$. With the condition $D$, $X$ is diffeomorphic to $\R^{+}\setminus\{0\}\times\S^{n-1}$ and $SO(n)$ acts on the factor $\S^{n-1}$. Taking, on $\S^{n-1}$, the hyperspherical coordinates $(\phi_{1},\ldots,\phi_{n-1})\in[-\frac{\pi}{2},\frac{\pi}{2}]^{n-2}\times[0,2\pi)$, we can consider the action of $\S^{1}$ on $X$, corresponding to one of the inclusions $\S^1\subset SO(n)$, defined by 
\[e^{i\theta}.(r,\phi_{1}\ldots,\phi_{n-1})=(r,\phi_{1}\ldots,\phi_{n-1}-\theta).\]  
If we consider the components $V_{k}\in L^{2}(\R^{+})\otimes H^{k}_{\omega_{max}^{k}}$ of $V$, they have the following form  \[V_{k}(r,\phi_{1}\ldots,\phi_{n-1})=s_{k}(r)v_{max}^{k}(\phi_{1}\ldots,\phi_{n-1})=s_{k}(r)(\prod_{i=1}^{n-2}\cos\phi_{i})^{k}e^{ik\phi_{n-1}}.\]
In fact we have $v_{max}^{k}=(x_{1}+ix_{2})^{k}$ with $x_{1}=(\underset{i=1}{\overset{n-2}{\prod}}\cos\phi_{i})\cos\phi_{n-1}$ and $x_{2}=(\underset{i=1}{\overset{n-2}{\prod}}\cos\phi_{i})\sin\phi_{n-1}$. So $V_{k}$ is $\S^{1}$ homogeneous of weight $k$ for the previously described $\S^{1}$ action.\\
\indent In conclusion the family of potentials constructed thanks to the action of $SO(n)$ is included into the potentials constructed with the action of $\S^{1}$.

So, why look at the $SO(n)$ action ? In fact, as we will see, using the $SO(n)$ action simplifies the proof of isoresonance. In particular we don't need the lower bound of the spectrum of the Laplacian on functions of weight $j$ i.e. the proposition \ref{prop:spectre}. This allows us to add to the free Laplacian a real $SO(n)$-invariant potential $V_{0}$ not compactly supported but just decreasing at infinity in order to continue $R_{V_{0}}$ (compare with remark \ref{remV0}). 
\end{rem}

In the way to prove the theorem \ref{thmSOn}, first, note that $V_{k}$ maps $L^{2}(\R^{+})\otimes H^{\ell}_{\omega_{max}^{\ell}}$ into $L^{2}(\R^{+})\otimes H^{\ell+k}_{\omega_{max}^{\ell+k}}$. As for the action of $\S^{1}$, this shift will be the key of the proof.

\subsection{From $L^{2}(X)$ to $L^{2}(X)^{+}$}

Let $\chi\in C^{\infty}_{c}(X)$ invariant under the action of $SO(n)$. As for the circular action we begin studying  $\Res(\Delta +\chi V)$, on $D_{N}^{+}$.
We can write the Lipmann-Schwinger equation and get that if $\lambda_{0}\in\Res(\Delta+\chi V)\cap(D_{N}^{+}\setminus\Res(\Delta))$ then there exists a nontrivial $u\in
L^{2}(X)$ such that 
\[\big(I+\rho^{-N}\chi V\widetilde{R}_{0}(\lambda_{0})\rho^{N}\big)u=0.\]

We want to prove that we can choose $u$ in $L^{2}(X)^{+}$ :
\begin{lemme}
\label{lemmeL2+}
For $\lambda_{0}\in\Res(\Delta+\chi V)\cap(D_{N}^{+}\setminus\Res(\Delta))$, there exists a non trivial $w\in L^{2}(X)^{+}$ such
that \[\big(I+\rho^{-N}\chi V\widetilde{R}_{0}(\lambda_{0})\rho^{N}\big)w=0.\]
\end{lemme}

Proof : as $\rho^{-N}\chi V\widetilde{R}_{0}(\lambda_{0})\rho^{N}: L^{2}(X)\rightarrow L^{2}(X)$ is a compact operator, we have
$\mathcal{H}_{-1}:=Ker\big(I+\rho^{-N}\chi V\widetilde{R}_{0}(\lambda_{0})\rho^{N}\big)$ is finite dimensional. 

In addition, for all
$\xi\in\mathfrak{g}_{+}$, $D_{\xi}$ maps $\mathcal{H}_{-1}$ into itself. Indeed, by definition of highest weight vector, we have 
$D_{\xi}V=\underset{k=1}{\overset{\infty}{\sum}}D_{\xi}V_{k}=0$. Moreover, as the action of $SO(n)$ is isometric, $D_{\xi}$ commutes with the Laplacian and so with $\widetilde{R}_{0}(\lambda_{0})$. $D_{\xi}$ commutes also with $\rho$ and $\chi$ because they are $SO(n)$-invariant. So, if $u\in\mathcal{H}_{-1}$ then $u=-\rho^{-N}\chi V\widetilde{R}_{0}(\lambda_{0})\rho^{N}u$ and
\begin{align*}
D_{\xi}u&=-\rho^{-N}\chi D_{\xi}(V\widetilde{R}_{0}(\lambda_{0})\rho^{N}u)\\
&=-\rho^{-N}\chi (D_{\xi}(V)\widetilde{R}_{0}(\lambda_{0})\rho^{N}u+VD_{\xi}(\widetilde{R}_{0}(\lambda_{0})\rho^{N}u))\\
&=-\rho^{-N}\chi V\widetilde{R}_{0}(\lambda_{0})\rho^{N}D_{\xi}u,
\end{align*}
and finally $D_{\xi}u\in\mathcal{H}_{-1}$.\\
\indent So $\mathcal{H}_{-1}$ is a finite representation of $\mathfrak{g}_{+}$. Moreover $\mathfrak{g}_{+}$ is a nilpotent algebra. It's coming from the fact that there is only a finite number of positive roots of $\mathfrak{g}$ and from the following calculation : if $\xi\in\mathfrak{g}_{\alpha}$ and $\zeta\in\mathfrak{g}_{\beta}$ then $\ad(\xi)(\zeta)\in\mathfrak{g}_{\alpha+\beta}$. To see this, for all $\sigma\in\mathfrak{h}$, we have
\begin{align*}
\ad(\sigma)([\xi,\zeta])=[\sigma,[\xi,\zeta]]&=[\xi,[\sigma,\zeta]]+[[\sigma,\xi],\zeta]\\
&=[\xi,\beta(\sigma)\zeta]+[\alpha(\sigma)\xi,\zeta]\\
&=(\alpha+\beta)(\sigma)[\xi,\zeta].
\end{align*}
Then by Engel's theorem (see \cite{FH} p.125) there exists a nonzero vector $w\in\mathcal{H}_{-1}$ such that $D_{\xi}w=0$ for all $\xi\in\mathfrak{g}_{+}$.\\
\indent We can decompose $w$ :
\[w=\sum_{k\in\N}w_{k},   \;\;w_{k}\in L^{2}(\R^{+})\otimes
H^{k},\]
and for $\xi\in\mathfrak{g}_{+}$ we have
\[D_{\xi}w=\sum_{k\in\N}D_{\xi}w_{k}=0,\]
with $D_{\xi}w_{k}\in L^{2}(\R^{+})\otimes H^{k}$. As the previous sum is direct, we have, for all $k$, $w_{k}\in (L^{2}(\R^{+})\otimes H^{k})\bigcap(\underset{\xi\in\mathfrak{g}_{+}}{\cap} \Ker D_{\xi})$. But, with the lemma \ref{lemmeGu}, we have  $(L^{2}(\R^{+})\otimes H^{k})\bigcap(\underset{\xi\in\mathfrak{g}_{+}}{\cap} \Ker D_{\xi})=L^{2}(\R^{+})\otimes H^{k}_{\omega_{max}^{k}}$. In conclusion, for all $k$, $w_{k}\in H^{k}_{\omega_{max}^{k}}$ and $w\in L^{2}(X)^{+}$. \ \ $\square$

\subsection{The end of the proof of theorem \ref{thmSOn}}

We assume that there exists $\lambda_{0}\in\Res(\Delta+\chi V)\cap(D_{N}^{+}\setminus\Res(\Delta))$, and we take a non trivial $w\in L^{2}(X)^{+}$ whose existence is given by the lemma \ref{lemmeL2+} and which satisfies
\[\big(I+\rho^{-N}\chi V\widetilde{R}_{0}(\lambda_{0})\rho^{N}\big)w=0.\]
For all $j\in\N$, we denote $w_{j}:=P_{j}w\in L^{2}(\R^{+})\otimes H^{j}_{\omega_{max}^{j}}$. The $P_{j}$ commute with $\rho$, $\chi$ because they are both invariant under the action of $SO(n)$ and $\widetilde{R}_{0}(\lambda_{0})$ because the action is isometric. We have
\[w_{j}= P_{j}\big(- \rho^{-N}\chi V\widetilde{R}_{0}(\lambda_{0})\rho^{N}w\big) = -\sum_{k=1}^{\infty}\rho^{-2N}\chi P_{j}\big(V_{k} \rho^{N}\widetilde{R}_{0}(\lambda_{0})\rho^{N}w\big).\]
Moreover we have already seen that, for all $\xi\in\mathfrak{g}_{+}$, $D_{\xi}$ commute with $\widetilde{R}_{0}(\lambda_{0})$ and $\rho$. Thus, if $w\in L^{2}(X)^{+}$ then $\rho^{N}\widetilde{R}_{0}(\lambda_{0})\rho^{N}w\in L^{2}(X)^{+}$ and we have $\rho^{N}\widetilde{R}_{0}(\lambda_{0})\rho^{N}w=\underset{\ell=0}{\overset{\infty}{\sum}}P_{\ell}(\rho^{N}\widetilde{R}_{0}(\lambda_{0})\rho^{N}w)$.

We use the shift created by highest weight vectors i.e. :
\[V_{k} : L^{2}(\R^{+})\otimes H^{\ell}_{\omega_{max}^{\ell}}\rightarrow L^{2}(\R^{+})\otimes H^{\ell+k}_{\omega_{max}^{\ell+k}}.\]
So
\[
w_{j}=-\sum_{k=1}^{\infty}\rho^{-2N}\chi V_{k} P_{j-k}\big(\rho^{N}\widetilde{R}_{0}(\lambda_{0})\rho^{N}w\big)= -\sum_{k=1}^{\infty}\rho^{-2N}\chi V_{k}\rho^{N}\widetilde{R}_{0}(\lambda_{0})\rho^{N}P_{j-k}(w).
\]
Thanks to the hypothesis $\underset{k}{\sup}\parallel\!V_{k}\!\parallel_{\infty}<+\infty$, the operators $\rho^{-2N}\chi V_{k}\rho^{N}\widetilde{R}_{0}(\lambda_{0})\rho^{N}$ are uniformly bounded in $k$ and consequently there exists a constant $C$ such that, for all $j\in\N$,
\[\parallel w_{j}\parallel\leq C\sum_{k=1}^{\infty}\parallel w_{j-k}\parallel.\] 
With this inequality and the fact that $w_{j}=0$ for all $j\leq 0$, we get $w_{j}=0$ for all $j\in\N$ and thus $w=0$ which is in contradiction with our hypothesis.

Finally we have proved that for all $\chi\in C^{\infty}_{c}(X)$ invariant under the action of $SO(n)$, we have $D_{N}^{+}\cap\Res(\Delta+\chi V)\subset\Res(\Delta)$.

In a second part we pass from $\chi V$ to $V$ like in the case of the $\S^{1}$ action. We introduce a family of smooth and compactly supported functions $(\chi_{r})$ invariant under the action of $SO(n)$ such that $\parallel\!(\chi_{r}-1)\rho\!\parallel_{\infty}$ tends to $0$ when $r$ tends to $+\infty$. We use the assumption $\rho^{-(N+1)}V\widetilde{R}_{0}(\lambda)\rho^{N}\in\mathcal{S}_{q}$ in order to characterize the resonances of $\Delta+V$ as the zeros of the holomorphic function in $\lambda$, $F(V,\lambda):=\det_{q}\big(I+\rho^{-N}V\widetilde{R}_{0}(\lambda)\rho^{N}\big)$. Finally with the Rouch\' e's theorem we prove that if $\widetilde{R}_{V}$ has a pole in $D_{N}^{+}\setminus\Res(\Delta)$ then, there exists $r$ such that $\widetilde{R}_{\chi_{r}V}$ has also a pole which is in contradiction with the previous part. So $D_{N}^{+}\cap\Res(\Delta+V)\subset\Res(\Delta)$.

To conclude, we prove $D_{N}^{+}\cap\Res(\Delta)\subset\Res(\Delta+V)$ and, in fact the equality with multiplicity, using the Agmon's perturbation theory of resonances exactly in the same way as in the case $\S^{1}$. This achieves the proof of theorem \ref{thmSOn}.

\section{Isoresonant potentials on the catenoid}

We are going to construct isoresonant potentials on the catenoid. In this case we use complex scaling defined by Wunsch and Zworski in \cite{WZ} instead of the Agmon's theory.

\subsection{Statement of the result}

The catenoid is the surface $X$ diffeomorphic to the cylinder $\R\times \S^{1}$ with the metric $g=\d r^{2}+(r^{2}+a^{2})\d\alpha^{2}$ where $(r,e^{i\alpha})\in\R\times\S^{1}$ and $a\in\R\setminus\{0\}$. We take $x=\frac{1}{\mid r\mid}$ outside $\{r=0\}$ as the function defining the boundary at infinity of $X$, $\partial_{\infty}X$, which is two copies of $\S^{1}$. Near the boundary, we have $g=\frac{\d x^{2}}{x^{4}}+\frac{(1+a^{2}x^{2})\d\alpha^{2}}{x^{2}}$ : it's a \textit{scattering metric} in the Melrose's sens (\cite{MR}).

The catenoid is an example in the Wunsch and Zworski article \cite{WZ} so we can use their results. They proved that there exists $\theta_{0}>0$ such that the resolvent of the free Laplacian, $(\Delta-z)^{-1}$, has a finite-meromorphic continuation from $\{z\in\C\ ;\ \Im z<0\}$ to $\{z\in\C\ ;\ \arg z<2\theta_{0}\}$ with values in operators from $L^{2}_{c}(X)$ to $H^{2}_{loc}(X)$. We denote $\widetilde{R}_{0}$ this continuation. Like in the previous part we call resonances its poles and we denote their set $\Res(\Delta)$.

The group $\S^{1}$ acts isometrically on $X$ by its trivial action on the factor  $\S^{1}$ : $e^{i\beta}.(r,e^{i\alpha})=(r,e^{i(\alpha+\beta)})$. Using this action we are going to construct isoresonant potentials :

\begin{thm}
\label{Potcatenoide}
Let $X$ be the catenoid $(\R\times \S^{1},\d r^{2}+(r^{2}+a^{2})\d\alpha^{2})$ with $(r,e^{i\alpha})\in \R\times \S^{1}$ and $a\in\R\setminus\{0\}$. We take $x=\frac{1}{\mid r\mid}$ outside $\{r=0\}$ as the function defining the boundary at infinity of $X$. 
Let $V\in xL^{\infty}(X)$ defined by
\[V(r,e^{i\alpha})=\sum_{m=1}^{\infty}V_{m}(r)e^{im\alpha}, \;\; (r,e^{i\alpha})\in\R\times \S^{1},\]
where, for all $m$, $V_{m}\in L^{\infty}(\R)$. We assume that, in a neighborhood of $\partial_{\infty}X$, $V(x,e^{i\alpha})$ with all its partial sums have a analytic continuation in $U\times W$ where $U$ is an open set of $\C$ including $\{\zeta\in\C\ ;\ \mid\!\zeta\!\mid\leq 1,\ 0\leq\arg \zeta\leq\theta_{0}\}$ with $\theta_{0}>0$ and $W$ is a neighborhood of $\S^{1}$ in $\C$. We also assume that, in $X$ and in all compacts of $U\times W$, the partial sums of $V$ tend to $V$ in infinite norm.

Then the resolvent $(\Delta+V-z)^{-1}$ has a finite-meromorphic continuation from $\{z\in\C\ ;\ \Im z<0\}$ to $\{z\in\C\ ;\ \arg z< 2\theta_{0}\}$ with values in operators from $L^{2}_{c}(X)$ to $H^{2}_{loc}(X)$, and moreover, in this open set, $\Res(\Delta+V)=\Res(\Delta)$ with the same multiplicities.

\end{thm}

We give an example of isoresonant potential on the catenoid :
\[V(x,e^{i\alpha})=\frac{xe^{i\alpha}}{1-\rho e^{i\alpha}}=x\sum_{m=1}^{\infty}\rho^{m-1}e^{im\alpha},\]
with $0\leq\rho<1$ and $U=\C$, $W=\{\omega\in\C\ ;\  \mid\!\omega\!\mid<\rho^{-1}\}.$

\subsection{Complex scaling}

We will use twice the complex scaling : to continue the resolvent of $\Delta+V$ and to study perturbations of resonances. So we begin with the description of this construction. We will follow \cite{WZ} where the construction is done for the free Laplacian. It is also valid when we add a potential $V\in x L^{\infty}(X)$.

We begin with the construction of a family $(X_{\theta})_{0\leq\theta\leq\theta_{0}}$ of submanifolds of $\C\times\C$ with the $\theta_{0}$ of the statement of theorem \ref{Potcatenoide}. They will be totally real i.e., for all $p\in X_{\theta}$, $T_{p}X_{\theta}\cap iT_{p}X_{\theta}=\{0\}$, and of maximal dimension. We define them as follows.

Let $\epsilon>0$ and $(t_{0},t_{1})\in]0,1[^{2}$ with $t_{0}<t_{1}$, then there exists a smooth deformation of $[0,1)$ in $U$, denote it $\gamma_{\theta}(t),  t\in[0,1)$ satisfying the following properties :
\begin{eqnarray}
\label{defgamma}
&\gamma_{\theta}(t)=te^{i\theta}  \;\;\;\textrm{for  $0\leq t<t_{0}$}&\nonumber\\
&\gamma_{\theta}(t)\equiv t  \;\;\; \textrm{for  $t>t_{1}$}&\nonumber\\
&\arg \gamma_{\theta}(t)\geq 0&\\
&0\leq \arg \gamma_{\theta}(t)-\arg \gamma_{\theta}'(t)\leq\epsilon&\nonumber\\
&0\leq 2 \arg \gamma_{\theta}(t)-\arg \gamma_{\theta}'(t)\leq\theta+\epsilon.\nonumber&
\end{eqnarray}

Now we can define $X_{\theta}:=(\gamma_{\theta}\times\partial_{\infty}X)\cup (X\cap\{x\geq 1\})$.

On a neighborhood of $\partial_{\infty}X$, the metric $g$ has the form $\frac{\d x^{2}}{x^{4}}+\frac{h}{x^{2}}$ where $h=(1+a^{2}x^{2})\d\alpha^{2}$ continues holomorphically to $U\times W$. So consider $P^{V}:=\Delta+V$, which is first an operator on $X$. If $V$ verifies the hypothesis of theorem \ref{Potcatenoide}, then its coefficients continue holomorphically in $U\times W$. We denote $\widetilde{P}^{V}$ the differential operator coming from this continuation. Since $X_{\theta}$ is totally real and of maximal dimension, we can define without ambiguity (cf \cite{SZ}) the differential operator $P^{V}_{\theta}$ by
\[\forall u\in C^{\infty}(X_{\theta}), \;\;  P^{V}_{\theta}u=(\widetilde{P}^{V}\widetilde{u})_{\mid X_{\theta}}\]
where $\widetilde{u}$ is an almost analytic extension of $u$, that is
\[\widetilde{u}\in C^{\infty}(U\times W), \;\; \widetilde{u}_{\mid X_{\theta}}=u, \;\; \overline{\partial}\widetilde{u}_{\mid X_{\theta}}=\mathcal{O}(\d(.,X_{\theta})^{N}), \;\; \textrm{for all $N$}.\]

We have 

\begin{proposition}
\label{spectrePtheta}
With $V\in xL^{\infty}(X)$, for all $0\leq\theta\leq\theta_{0}$, $P^{V}_{\theta}$ has a discrete spectrum in $\C\setminus e^{2i\theta}\R^{+}$. Moreover, for $\theta$ such that $0\leq\theta_{2}\leq\theta\leq\theta_{0}$, the spectrum of $P^{V}_{\theta}$ in $\{0\leq\arg z<2\theta_{2}\}$ with its multiplicity do not depend on $\theta$. This spectrum doesn't depend too on the choice of a $\gamma_{\theta}$ verifying (\ref{defgamma}).
\end{proposition}
Proof : Following exactly the Wunsch and Zworski proof in \cite{WZ} we can prove that, for all $z\in\C\setminus e^{2i\theta}\R^{+}$, $P^{V}_{\theta}-z : H^{2}(X_{\theta})\rightarrow L^{2}(X_{\theta})$ is a Fredholm operator with index zero. The unique difference is the presence of our potential $V$. But, since it is null at the boundary of $X_{\theta}$ ($=\partial_{\infty}X$), it doesn't change the principal symbol and the normal symbol of $\Delta_{\theta}-z$.

\subsection{Continuation of the resolvent}

We want to get the meromorphic continuation of the resolvent  $R_{V}(z):=(\Delta+V-z)^{-1}$ from $\{z\in\C\ ;\ \Im z<0\}$ to $\{z\in\C\ ;\ \arg z< 2\theta_{0}\}$ with values in operators from $L^{2}_{c}(X)$ to $H^{2}_{loc}(X)$.

Let $z$ with $\arg z<2\theta_{0}$ which is not an eigenvalue of $P^{V}_{\theta_{0}}$. Take $f\in L^{2}_{c}(X)$. With the proposition \ref{spectrePtheta}, we can choose $\gamma_{\theta_{0}}$ and more precisely the $t_{1}$ in the definition \ref{defgamma} such that, on the support of $f$, $X_{\theta_{0}}$ coincides with $X$. Then $f\in L^{2}(X_{\theta_{0}})$ and there exists an unique solution $u_{\theta_{0}}\in H^{2}(X_{\theta_{0}})$ of
\[(P^{V}_{\theta_{0}}-z)u_{\theta_{0}}=f.\]

We give a lemma whose proof is given in \cite{SZ},

\begin{lemme}
\label{lemmeprolong}
Let $\Omega\subset\C^{n}$ an open set, a compact $K\subset\Omega$ and a continuous family $X_{t}, t\in[0,1]$ of totally real submanifolds of $\Omega$ of maximal dimension such that $X_{t}\cap(\Omega\setminus K)=X_{t'}\cap(\Omega\setminus K)$ for all $t, t'\in [0,1]$. Let $\widetilde{P}$ a differential operator with holomorphic coefficients in $\Omega$ such that $P_{X_{t}}$ (the restriction of $\widetilde{P}$ on $X_{t}$) is elliptic for all $t\in[0,1]$. If $u$ is a distribution on $X_{0}$ and if $P_{X_{0}}u$ continues as an holomorphic function on a neighborhood of $\underset{t\in[0,1]}{\bigcup}X_{t}$ then the same is true for $u$. 
\end{lemme}

$f$ has an holomorphic continuation to $\underset{\theta\in[0,\theta_{0}]}{\bigcup}X_{\theta}$, because deformations occur outside its support. Since $P^{V}_{\theta}-z$ is elliptic for all $\theta\in[0,\theta_{0}]$, we can apply the previous lemma \ref{lemmeprolong}, and get an holomorphic continuation of $u_{\theta_{0}}$ on $\underset{\theta\in[0,\theta_{0}]}{\bigcup}X_{\theta}$. We denote by $G$ this continuation.\\
\indent Then we define the continuation to the resolvent by 
\[\widetilde{R}_{V}(z)f=G_{\mid\!X_{0}}\in H^{2}(X).\]

Now take $z_{0}$ an eigenvalue of $P^{V}_{\theta_{0}}$, then it is also an eigenvalue $P^{V}_{\theta}$ for all $\arg z_{0}<2\theta\leq 2\theta_{0}$. For $z$ near $z_{0}$ and $\theta$ such that $\arg z<2\theta$ we have the following Laurent expansion
\[(P^{V}_{\theta}-z)^{-1}=\sum_{j=1}^{M(z_{0})}\frac{A_{j}^{\theta}(z_{0})}{(z-z_{0})^{j}}+H_{\theta}(z,z_{0}),\] 
where $A_{j}^{\theta}(z_{0})$ are finite rank operators and $H_{\theta}(z,z_{0})$ is holomorphic in $z$ near $z_{0}$. Still following \cite{WZ}, we obtain that the continued resolvent has, near each of its pole, a Laurent expansion with exactly the same form.

In conclusion, a resonance $z_{0}\in\{\arg z< 2\theta_{0}\}$ of $\Delta_{X}+V$, which is first defined as poles of the continuation of the resolvent, is also characterized as an element of the spectrum of a $P^{V}_{\theta}$ with $\arg z_{0}<2\theta\leq 2\theta_{0}$. Multiplicities and orders are the same in the two visions, so, thanks to the proposition \ref{spectrePtheta}, they do not depend on the chosen $\theta$.

\subsection{Proof of the isoresonance}

\subsubsection{Localization of resonances for the truncated partial sums of $V$}

Let $S_{M}(r,e^{i\alpha})=\underset{m=1}{\overset{M}{\sum}}V_{m}(r)e^{im\alpha}$ and $\chi\in C^{\infty}_{c}(X)$, $\S^{1}$ invariant. In this part we will prove $\Res(\Delta +\chi S_{M})\subset\Res(\Delta)$ in $D^{+}:=\{z\in\C\ ;\ \arg z< 2\theta_{0}\}$.

We take another $\S^{1}$ invariant cutoff function  $\chi_{1}\in C^{\infty}_{c}(X)$ such that $\chi_{1}=1$ on the support of $\chi$. Then we have, for $z\in D^{+}\setminus\Res(\Delta)$,
\[(\Delta+\chi S_{M}-z)\widetilde{R}_{0}(z)\chi_{1}=\chi_{1}(I+\chi S_{M}\widetilde{R}_{0}(z)\chi_{1}).\]

$\chi S_{M}\widetilde{R}_{0}(z)\chi_{1}$ is an holomorphic family of compact operators in $D^{+}\setminus\Res(\Delta)$ such that 
\[\parallel\!\chi S_{M}\widetilde{R}_{0}(z)\chi_{1}\!\parallel<1\]
with $\mid\!z\!\mid$ sufficiently large in $\{z\in\C\ ;\ \Im z<0\}$.
So we can apply the Fredholm analytic theory and get $(I+\chi S_{M}\widetilde{R}_{0}(z)\chi_{1})^{-1}$ and thus $\widetilde{R}_{\chi S_{M}}(z):=(\Delta+\chi S_{M}-z)^{-1}$ meromorphic in $D^{+}\setminus\Res(\Delta)$. Moreover, in $D^{+}\setminus\Res(\Delta)$, we can characterize poles of $\widetilde{R}_{\chi S_{M}}$, that is resonances, with the existence of a non trivial $u\in L^{2}(X)$ solution of
\[(I+\chi S_{M}\widetilde{R}_{0}(z)\chi_{1})u=0.\]

After, we prove that this non trivial solution $u$ can't exist. It is exactly the same proof than for theorem \ref{thm : circular}. We use the shift created by the components $V_{m}(r)e^{im\alpha}$ on the spaces $L^{2}_{j}(X)$. There is just two things to verify. First the catenoid satisfies the condition $C$. In fact, all compact $K$ of $X$ is included in a compact manifold with boundary $\widetilde{K}=[-R,R]\times\S^{1}$ on which we can put the metric $\widetilde{g}=\d r^{2}+f(r)\d\alpha^{2}$ with $f(r)=r^{2}+a^{2}$ in $K$ and constant near the boundary of $\widetilde{K}$. We also have to verify that $\widetilde{R}_{0}$ commutes with the action of $\S^{1}$ in order to have the commutation with the projectors $P_{j}$. For that, remember that the complex scaling doesn't touch the factor $\partial_{\infty}X$ of the catenoid and $\S^{1}$ only acts on this factor. So the action of $\S^{1}$ is isometric on all the $X_{\theta}$ and so it commutes with $\widetilde{R}_{0}$.

Finally we get $\Res(\Delta +\chi S_{M})\subset\Res(\Delta)$ in $D^{+}$.

\subsubsection{Localization of resonances for $V$}

We have to control perturbations of resonances when we pass from $\chi S_{M}$ to $V$. Instead of use regularized determinants as before (they were adapted to the weight spaces $\rho^{N} L^{2}$ but not to cutoff functions), we use complex scaling in order to transform resonances into eigenvalues.

Assume, ab absurdum, $\Delta+V$ has a resonance $z_{0}$ in $D^{+}\setminus\Res(\Delta)$. Using complex scaling, it means : $z_{0}$ is an eigenvalue of $P^{V}_{\theta}$ with $\arg z_{0}<2\theta\leq 2\theta_{0}$. Let $\Omega\subset\{z\in\C\ ;\ \arg z<2\theta\}\setminus~\Res(\Delta)$ an open set, with a smooth boundary $\Gamma$, containing $z_{0}$ and such that $\overline{\Omega}\cap\Res(\Delta+V)=\{z_{0}\}$. Our aim is to show that there exist a $\S^{1}$ invariant cutoff function $\chi$, and $M$ such that $P^{\chi S_{M}}_{\theta}$ has an eigenvalue in $\Omega$. If we do that, $\Delta +\chi S_{M}$ would have a resonance in $\Omega$ which would be in contradiction with the previous part.

We have assumed in theorem \ref{Potcatenoide} that, for all $M$, $S_{M}$ has an analytic continuation to $U\times W$, and $V$ too. Hence we can restrict these two continuations to $X_{\theta}$ and now work on $X_{\theta}$. $V\in x L^{\infty}(X)$, so $V$ tends to $0$ when we reach $\partial X_{\theta}$. Consequently there exists $\chi\in C^{\infty}_{c}(X_{\theta})$, $\S^{1}$ invariant, which continues analytically in $U\times W$, such that $\parallel\!\chi V-V\!\parallel_{L^{\infty}(X_{\theta})}$ is as small as we want. With the hypothesis of theorem \ref{Potcatenoide}, we also have that the partial sums  $S_{M}$ tend to $V$ on $X_{\theta}$ in infinite norm. Finally there exist $\chi$ like we have just described, and $M$ such that 
\[\parallel\!V-\chi S_{M}\!\parallel_{L^{\infty}(X_{\theta})}<\frac{\delta^{2}}{\delta+\frac{\ell}{2\pi}}\]
where $\delta^{-1}=\underset{z\in\Gamma}{\max}\parallel\!(P^{V}_{\theta}-z)^{-1}\!\parallel$ and $\ell$ is the length of $\Gamma$. We have been inspired by Gohberg and Krejn in \cite{GK} (theorem 3.1). 

We consider the projectors of $L^{2}(X_{\theta})$ associated with the generalized eigenspaces of the two operators that we are comparing :
\[\begin{array}{rcl}
\Pi_{V}&=&\frac{1}{2\pi i}\int_{\Gamma}(P^{V}_{\theta}-z)^{-1}\d z\\
\Pi_{\chi S_{M}}&=&\frac{1}{2\pi i}\int_{\Gamma}(P^{\chi S_{M}}_{\theta}-z)^{-1}\d z
\end{array}\]
We have $(P^{\chi S_{M}}_{\theta}-z)^{-1}=(P^{V}_{\theta}-z)^{-1}\big(I+(\chi S_{M}-V)(P^{V}_{\theta}-z)^{-1}\big)^{-1}$. Since $\frac{\delta^{2}}{\delta+\frac{\ell}{2\pi}}<\delta$, for all $z\in\Gamma$, we can be certain of the convergence in
\[(P^{\chi S_{M}}_{\theta}-z)^{-1}=(P^{V}_{\theta}-z)^{-1}\big(I+\sum_{j=1}^{\infty}[(V-\chi S_{M})(P^{V}_{\theta}-z)^{-1}]^{j}\big).\]
Look at the difference between the two projectors :
\[\Pi_{\chi S_{M}}-\Pi_{V}=\frac{1}{2\pi i}\int_{\Gamma}(P^{V}_{\theta}-z)^{-1}\sum_{j=1}^{\infty}[(V-\chi S_{M})(P^{V}_{\theta}-z)^{-1}]^{j}\d z.\]
Hence
\[\parallel\!\Pi_{\chi S_{M}}-\Pi_{V}\!\parallel\leq \frac{\ell}{2\pi}\max_{z\in\Gamma}\frac{\parallel\!(P^{V}_{\theta}-z)^{-1}\!\parallel^{2} \parallel\!V-\chi S_{M}\!\parallel_{L^{\infty}(X_{\theta})}}{1-\parallel\!(P^{V}_{\theta}-z)^{-1}\!\parallel \parallel\!V-\chi S_{M}\!\parallel_{L^{\infty}(X_{\theta})}},\]
but $\parallel\!(P^{V}_{\theta}-z)^{-1}\!\parallel\leq\delta^{-1}$ by definition of $\delta$ and $\parallel\!V-\chi S_{M}\!\parallel_{L^{\infty}(X_{\theta})}<\frac{\delta^{2}}{\delta+\frac{\ell}{2\pi}}$ by hypothesis, hence $1-\parallel\!(P^{V}_{\theta}-z)^{-1}\!\parallel \parallel\!V-\chi S_{M}\!\parallel_{L^{\infty}(X_{\theta})}>1-\frac{\delta}{\delta+\frac{\ell}{2\pi}}$ and so
\[\parallel\!\Pi_{\chi S_{M}}-\Pi_{V}\!\parallel<1.\]

Consequently the ranges of $\Pi_{V}$ and $\Pi_{\chi S_{M}}$ have the same dimension, which is not zero because $z_{0}$ is an eigenvalue of $P_{\theta}^{V}$, so $P^{\chi S_{M}}_{\theta}$ has an eigenvalue in $\Omega$ and we have our contradiction.

\subsubsection{Persistence of resonances}

To finish the proof of theorem \ref{Potcatenoide} we have to show that $\Res(\Delta)\subset\Res(\Delta+V)$ in $D^{+}=\{z\in\C\ ;\ \arg z<2\theta_{0}\}$. This time we use the complex scaling instead of the Agmon's perturbation theory of resonances.

Let $z_{0}\in D^{+}$ a resonance of $\Delta$ with multiplicity $m$. We take $\Omega\subset\{z\in\C\ ;\ \arg z<2\theta_{0}\}$ a domain such that $\overline{\Omega}\cap\Res(\Delta)=\{z_{0}\}$. We consider the family of operators $\Delta+tV$ with $t\geq 0$. Remark that $tV$ verifies the hypothesis of theorem \ref{Potcatenoide}, so we can localize its resonances as in the previous part, $\Res(\Delta+tV)\subset\Res(\Delta)$, and thus :
\[\Res(\Delta+tV)\cap\Omega\subset\{z_{0}\}.\] 

We are going to prove, by connexity, that the following set 
\[E:=\{t_{0}\geq 0\ ;\ \forall t\in[0,t_{0}],\ \Res(\Delta+tV)\cap\Omega=\{z_{0}\} \ \ \textrm{with multiplicity}\  m\},\]
is equal to $[0,+\infty[$. It is not empty because $0\in E$.

Take $t_{0}\in E$, and $\theta$ such that $\arg z_{0}<2\theta\leq 2\theta_{0}$ and $\Omega\subset\{z\in\C\ ;\ \arg z< 2\theta\}$. We know by complex scaling that the spectrum of $P_{\theta}^{t_{0}V}$ in $\Omega$ exactly corresponds with the resonances of $\Delta+t_{0}V$ with the same multiplicities. Hence 
\[\spec(P_{\theta}^{t_{0}V})\cap\Omega=\{z_{0}\}.\] 
Moreover $t\rightarrow P_{\theta}^{tV}$ is an holomorphic family in the sens of Kato for $t$ in a complex neighbourhood of $t_{0}$ because $P_{\theta}^{tV}=\Delta_{\theta}+tV_{\mid X_{\theta}}$ and $V$ is bounded in $X_{\theta}$. So its eigenvalues are continuous for $t$ in a neighbourhood of $t_{0}$. But with the localization of the resonances of $\Delta+tV$ we obtain that for all $t$
\[\spec(P_{\theta}^{tV})\cap\Omega\subset\{z_{0}\},\]
so, there exists $\varepsilon>0$, such that for all $t\in]t_{0}-\varepsilon,t_{0}+\varepsilon[$,  
\[\spec(P_{\theta}^{tV})\cap\Omega=\{z_{0}\},\]
with multiplicity $m$. Then, thanks to the complex scaling parallel, we get that, for all $t\in]t_{0}-\varepsilon,t_{0}+\varepsilon[$
\[\Res(\Delta+tV)\cap\Omega=\{z_{0}\},\]
with multiplicity $m$ and thus $]t_{0}-\varepsilon,t_{0}+\varepsilon[$ is included in $E$ which is open.\\
\indent We show that $E$ is closed too doing the same work with the complementary set of $E$ in $[0,+\infty[$.

In conclusion $E=[0,+\infty[$ and taking $t=1$ we have, in $\Omega$, $\Res(\Delta+V)=\Res(\Delta)$ with the same multiplicities. Doing the same work in the neighbourhood of each resonance of the free Laplacian we obtain $\Res(\Delta+V)=\Res(\Delta)$ with the same multiplicities in all $D^{+}$ which finishes the proof of theorem \ref{Potcatenoide}.

\bibliography{bibisoART}

\end{document}